%% file: main.tex
\documentclass[authoryear,a4paper,11pt,preprint,3p]{elsarticle}
\usepackage[onehalfspacing]{setspace}
\usepackage{wrapfig}
\usepackage{graphicx}
\usepackage{tabularx}
\usepackage{multirow}
\usepackage{mathtools}
\usepackage{amssymb}
\usepackage{bm}
\usepackage{cases}
\usepackage{caption}
\usepackage{subcaption}
\usepackage{threeparttable}
\usepackage[x11names]{xcolor}
\usepackage{transparent}
\usepackage{soul}
\usepackage{enumitem}
\usepackage{array}
\usepackage{listings}
\usepackage{colortbl}
\usepackage{adjustbox}
\usepackage{natbib}
\usepackage{har2nat}
\usepackage{hyperref}
\usepackage{enumitem}
\usepackage{framed}
\usepackage{pdflscape}
\usepackage{csvsimple-l3}
\usepackage{longtable}
\usepackage{acronym}
\usepackage[byauthor]{credits}
\usepackage[symbol]{footmisc}

\hypersetup{colorlinks=true}
\hypersetup{linkcolor=blue}
\makeatletter
\renewcommand\hyper@natlinkbreak[2]{#1}
\makeatother

\newcommand{\legendColor}{gray!10}
\newcommand{\tableHeaderColor}{gray!10}
\newcommand{\matrixBlockColor}{blue!10}
\newcommand{\correspondingauthor}{\footnote{Corresponding Author: l.hadidi@fz-juelich.de}}
\newcolumntype{R}[1]{>{\raggedright\arraybackslash}m{#1}}
\hypersetup{pdfsubject={Systematic Review on Parallelization Strategies via Decomposition}}
\hypersetup{pdfkeywords={OR in energy, Large scale optimization, Combinatorial optimization, Linear programming, High Performance Computing}}
\journal{arXiv}
\begin{document}
    \begin{frontmatter}
    
        \title{Large-Scale Linear Energy System Optimization: A Systematic Review on Parallelization Strategies via Decomposition}
        
        \author[a]{Lars Hadidi \correspondingauthor}
        \author[c]{Leonard Göke}
        \author[a]{Maximilian Hoffmann}
        \author[d]{Mario Klostermeier}
        \author[f]{Shima Sasanpour}
        \author[e]{Tim Varelmann}
        \author[d]{Vassilios Yfantis}      
        \author[a]{Jochen Linßen}
        \author[a,b]{Detlef Stolten}
        \author[a]{Jann M. Weinand}
        
        \affiliation[a]{organization={Forschungszentrum Jülich GmbH, Institute of Energy and Climate Research – Jülich Systems Analysis (ICE-2), 52425 Jülich, Germany}}
        
        \affiliation[b]{organization={RWTH Aachen University, Chair for Fuel Cells, Faculty of Mechanical Engineering, 52062 Aachen, Germany}}
        
        \affiliation[c]{organization={ETH Zurich, Reliability and Risk Engineering, Department of Mechanical and Process Engineering, 8092 Zurich, Switzerland}}

        \affiliation[d]{organization={RPTU Kaiserslautern-Landau, Chair of Machine Tools and Control Systems, Department of Mechanical and Process Engineering, 67663 Kaiserslautern, Germany}}
        
        \affiliation[e]{organization={Bluebird Optimization, 48429 Rheine, Germany}}

        \affiliation[f]{organization={German Aerospace Center (DLR), Institute of Networked Energy Systems, 70563 Stuttgart, Germany}}

        \credit{Lars Hadidi}
        {Conceptualization,Methodology,Validation,Formal analysis,Investigation,Data curation,Writing -- original draft,Writing -- review \& editing,Visualization}
        
        \credit{Leonard Göke}
        {Conceptualization,Investigation,Data curation,Validation}
        
        \credit{Maximilian Hoffmann}
        {Conceptualization,Investigation,Resources}
        
        \credit{Mario Klostermeier}
        {Writing -- original draft}
        
        \credit{Shima Sasanpour}
        {Writing -- original draft,Validation}
        
        \credit{Tim Varelmann}
        {Investigation,Data curation,Software,Writing -- original draft}
        
        \credit{Vassilios Yfantis}
        {Conceptualization}  
        
        \credit{Jochen Linßen}
        {Resources,Funding acquisition}
        
        \credit{Detlef Stolten}
        {Resources,Funding acquisition}

        \credit{Jann Weinand}
        {Conceptualization,Methodology,Writing -- review \& editing,Supervision,Project administration}
        
        \begin{abstract}
            As renewable energy integration, sector coupling, and spatiotemporal detail increase, energy system optimization models grow in size and complexity, often pushing solvers to their performance limits. This systematic review explores parallelization strategies that can address these challenges. We first propose a classification scheme for linear energy system optimization models, covering their analytical focus, mathematical structure, and scope. We then review parallel decomposition methods, finding that while many offer performance benefits, no single approach is universally superior. The lack of standardized benchmark suites further complicates comparison. To address this, we recommend essential criteria for future benchmarks and minimum reporting standards. We also survey available software tools for parallel decomposition, including modular frameworks and algorithmic abstractions. Though centered on energy system models, our insights extend to the broader operations research field.             
        \end{abstract}
        
        
    \end{frontmatter}
    
    \input{sections/01-introduction}
    \input{sections/02-theory}
    \input{sections/03-methods}
    \input{sections/04-energySystemModels}
    \input{sections/05-parallelization}
    \input{sections/06-implementations}
    \input{sections/07-conclustion}

    \newpage
    \section*{CRediT Author Statement}
    \insertcreditsstatement 

    \section*{Abbreviations}\label{sec:acronyms}
    \small
    \begin{acronym}[CCCCC]
    \acro{BB}{Branch-And-Bound}
    \acro{BPC}{Branch-And-Price-And-Cut}
    \acro{DecSys}{Decomposed System}
    \acro{ESOM}{Energy System Optimization Model}
    \acro{LP}{Linear Program}
    \acro{MEnv}{Modelling Environment}
    \acro{MILP}{Mixed Integer Linear Program}
    \acro{MIP}{Mixed Integer Program}
    \acro{MPI}{Message Passing Interface}
    \acro{PMet}{Performance Metric}
    \acro{PVal}{Performance Value}
    \acro{QMet}{Quality Metric}
    \acro{QVal}{Quality Value}
    \acro{SIMD}{Single Instruction Multiple Data}
    \acro{SP}{Stochastic Program}    
    \end{acronym}
    \normalsize
    
    \newpage
    \tiny
    \bibliographystyle{bibliography/style}
    \bibliography{bibliography/bibliography}
    \normalsize

    \appendix
    \section{Literature Data}\label{app:data}
    \begin{landscape}
    \tiny
    \begin{longtable}{|
    m{.06\linewidth}|
    m{.08\linewidth}|
    m{.05\linewidth}|
    m{.05\linewidth}|
    m{.05\linewidth}|
    m{.06\linewidth}|
    m{.06\linewidth}|
    m{.06\linewidth}|
    m{.06\linewidth}|
    m{.06\linewidth}|
    m{.06\linewidth}|
    m{.07\linewidth}|
    m{.04\linewidth}|
    m{.04\linewidth}|}
        \caption{Method Benchmarks. \textit{Abbreviations: \acf{ESOM}, \acf{MEnv}, \acf{PMet}, \acf{QMet}, \acf{DecSys}, \acf{PVal}, \acf{QVal}}}
        \label{tab:results-method}
        \\ \hline
        \cellcolor{\tableHeaderColor}\textbf{Study} & 
        \cellcolor{\tableHeaderColor}\textbf{\acs{ESOM} Name} & 
        \cellcolor{\tableHeaderColor}\textbf{\acs{ESOM} Type} &
        \cellcolor{\tableHeaderColor}\textbf{\acs{ESOM} Size} &
        \cellcolor{\tableHeaderColor}\textbf{\acs{MEnv}} &
        \cellcolor{\tableHeaderColor}\textbf{Reference Solver} &
        \cellcolor{\tableHeaderColor}\textbf{\acs{PMet}} &
        \cellcolor{\tableHeaderColor}\textbf{\acs{QMet}} &
        \cellcolor{\tableHeaderColor}\textbf{Reference \acs{PVal}} &
        \cellcolor{\tableHeaderColor}\textbf{Reference \acs{QVal}} &
        \cellcolor{\tableHeaderColor}\textbf{\acs{DecSys} Method} &
        \cellcolor{\tableHeaderColor}\textbf{Decomposed Dimensions} &
        \cellcolor{\tableHeaderColor}\textbf{\acs{DecSys} \acs{PVal}} &
        \cellcolor{\tableHeaderColor}\textbf{\acs{DecSys} \acs{QVal}}
        \begin{filecontents*}{results/DataExtraction-MethodBenchmark.csv}
PublicationTitle,Citation,ESOMName,ESOMType,ESOMSize,AML,BaselineSolver,QualityMetric,PerformanceMetric,BaselinePerformance,BaselineQuality,DecompositionMethod,DecompositionDimensions,DecompositionPerformance,DecompositionQuality,SpeedUp,Continuity,Integrality,Stochasticity,SHDecisions,EXDecisions
A Consensus-ADMM Approach for Strategic Generation Investment in Electricity Markets,\cite{Dvorkin2018},Custom Strategic Investment,BU-CI-SHEX.t,{9,648 variables},GAMS,CPLEX 12,Custom Gap,Runtime in seconds,3624.0,0 %,Lagrange Consensus ADMM,scenarios,9.0,0.5 %,402.7,$\blacksquare$,$\blacksquare$,{},$\blacksquare$,$\blacksquare$
Nesterov-Type Accelerated ADMM (N-ADMM) with Adaptive Penalty for Three-Phase Distributed OPF under Non-Ideal Data Transfer Scenarios,\cite{Paul2023},IEEE 123-Bus,BU-CI-SH.st,-,MATLAB,MATLAB 2019,Relative Baseline Deviation,Runtime in seconds,51.0,0 %,Lagrange Distributed ADMM,spatial,193.0,0.208 %,0.3,$\blacksquare$,$\blacksquare$,{},$\blacksquare$,
Nesterov-Type Accelerated ADMM (N-ADMM) with Adaptive Penalty for Three-Phase Distributed OPF under Non-Ideal Data Transfer Scenarios,\cite{Paul2023},IEEE 123-Bus,BU-CI-SH.st,-,MATLAB,MATLAB 2019,Relative Baseline Deviation,Runtime in seconds,51.0,0 %,Lagrange Distributed Adaptive N-ADMM,spatial,54.0,0.136 %,0.9,$\blacksquare$,$\blacksquare$,{},$\blacksquare$,
Computational study of security constrained economic dispatch with multi-stage rescheduling,\cite{Liu2015},IEEE 2383-Bus WinterPeak 2349c,BU-C-SHEX.st,-,GAMS,CPLEX -,Convergence,Runtime in seconds,7200.0,no,Benders Fatmaster Heuristic,scenarios,516.0,-,14.0,$\blacksquare$,{},,$\blacksquare$,$\blacksquare$
Computational study of security constrained economic dispatch with multi-stage rescheduling,\cite{Liu2015},IEEE 2736-Bus SummerPeak 2749c,BU-C-SHEX.st,-,GAMS,CPLEX -,Convergence,Runtime in seconds,7200.0,no,Benders Fatmaster Heuristic,scenarios,221.0,-,32.6,$\blacksquare$,{},,$\blacksquare$,$\blacksquare$
Computational study of security constrained economic dispatch with multi-stage rescheduling,\cite{Liu2015},IEEE 2737-Bus SummerOffPeak 2753c,BU-C-SHEX.st,-,GAMS,CPLEX -,Convergence,Runtime in seconds,7200.0,no,Benders Fatmaster Heuristic,scenarios,101.0,-,71.3,$\blacksquare$,{},,$\blacksquare$,$\blacksquare$
Computational study of security constrained economic dispatch with multi-stage rescheduling,\cite{Liu2015},IEEE 2746-Bus WinterOffPeak 2794c,BU-C-SHEX.st,-,GAMS,CPLEX -,Convergence,Runtime in seconds,7200.0,no,Benders Fatmaster Heuristic,scenarios,119.0,-,60.5,$\blacksquare$,{},,$\blacksquare$,$\blacksquare$
Computational study of security constrained economic dispatch with multi-stage rescheduling,\cite{Liu2015},IEEE 2746-Bus WinterPeak 2719c,BU-C-SHEX.st,-,GAMS,CPLEX -,Convergence,Runtime in seconds,7200.0,no,Benders Fatmaster Heuristic,scenarios,334.0,-,21.6,$\blacksquare$,{},,$\blacksquare$,$\blacksquare$
Modeling and Solution of the Large-Scale Security-Constrained Unit Commitment,\cite{Fu2013},IEEE 1168-Bus 168h,BU-CI-SH.st,-,-,CPLEX 12,Optimality Gap,Runtime in seconds,10500.0,-,Lagrange,{operational, scenarios},350.0,-,30.0,$\blacksquare$,$\blacksquare$,{},$\blacksquare$,
A hybrid approach for unit commitment with splitting technique and local search,\cite{Zhang2024},Custom 1080 Generators UC 48p,BU-CI-SH.st,-,MATLAB,CPLEX 12,Optimality Gap,Runtime in seconds,7200.0,7.02 %,Lagrange ADMM,temporal,110.0,-,65.5,$\blacksquare$,$\blacksquare$,{},$\blacksquare$,
A hybrid approach for unit commitment with splitting technique and local search,\cite{Zhang2024},Custom 1080 Generators UC 96p,BU-CI-SH.st,-,MATLAB,CPLEX 12,Optimality Gap,Runtime in seconds,7200.0,7.47 %,Lagrange ADMM,temporal,219.0,-,32.9,$\blacksquare$,$\blacksquare$,{},$\blacksquare$,
A hybrid approach for unit commitment with splitting technique and local search,\cite{Zhang2024},Custom 1080 Generators UC 168p,BU-CI-SH.st,-,MATLAB,CPLEX 12,Optimality Gap,Runtime in seconds,7200.0,-,Lagrange ADMM,temporal,875.0,-,8.2,$\blacksquare$,$\blacksquare$,{},$\blacksquare$,
Accelerating the benders decomposition for network-constrained unit commitment problems,\cite{Wu2010},Custom 5663-Bus NCUC Peak-Hour,BU-CI-SH.st,{21,115 variables},-,CPLEX 11,Convergence,Runtime in seconds,67.0,yes,Benders Strong Multi-Cut,{operational, spatial},25.0,yes,2.7,$\blacksquare$,$\blacksquare$,{},$\blacksquare$,
Accelerating the benders decomposition for network-constrained unit commitment problems,\cite{Wu2010},Custom 5663-Bus NCUC 24-Hours,BU-CI-SH.st,{500,736 variables},-,CPLEX 11,Convergence,Runtime in seconds,10800.0,no,Benders Strong Multi-Cut,{operational, spatial},2159.0,yes,5.0,$\blacksquare$,$\blacksquare$,{},$\blacksquare$,
A comparison of implicit and explicit methods for contingency constrained unit commitment,\cite{Huang2017},IEEE 24-Bus,BU-CI-SH.st,{1,889,569 variables},AMPL,CPLEX 12,Convergence,Runtime in seconds,259200.0,no,Benders Type Explicit Constraint Sets,operational,228.7,yes,1133.4,$\blacksquare$,$\blacksquare$,{},$\blacksquare$,
Solving Combined Sizing and Dispatch of PV and Battery Storage for a Microgrid using ADMM,\cite{Steven2024},Custom LV Microgrid,BU-CI-SHEX.st,{338,400 variables},JuMP,Gurobi 10,Relative Baseline Deviation,Runtime in seconds,14.5,0 %,Lagrange ADMM,temporal,3.5,0 %,4.1,$\blacksquare$,$\blacksquare$,{},$\blacksquare$,$\blacksquare$
Decomposing a renewable energy design and dispatch model,\cite{Wales2024},REopt (CHP),BU-CI-SHEX.te,{227,800 to 403,100 variables},AMPL,CPLEX 12,Optimality Gap,Runtime in seconds,300.0,13.9 %,Lagrange,temporal,300.0,7.1 %,1.0,$\blacksquare$,$\blacksquare$,{},$\blacksquare$,$\blacksquare$
Decomposing a renewable energy design and dispatch model,\cite{Wales2024},{REopt (CHP, TES)},BU-CI-SHEX.te,{227,800 to 403,100 variables},AMPL,CPLEX 12,Optimality Gap,Runtime in seconds,300.0,10.7 %,Lagrange,temporal,300.0,8.3 %,1.0,$\blacksquare$,$\blacksquare$,{},$\blacksquare$,$\blacksquare$
Decomposing a renewable energy design and dispatch model,\cite{Wales2024},{REopt (CHP, PV, BES)-5b},BU-CI-SHEX.te,{227,800 to 403,100 variables},AMPL,CPLEX 12,Optimality Gap,Runtime in seconds,300.0,83.7 %,Lagrange,temporal,300.0,2.0 %,1.0,$\blacksquare$,$\blacksquare$,{},$\blacksquare$,$\blacksquare$
Decomposing a renewable energy design and dispatch model,\cite{Wales2024},{REopt (CHP, PV, CHILL, BES, TES)},BU-CI-SHEX.te,{227,800 to 403,100 variables},AMPL,CPLEX 12,Optimality Gap,Runtime in seconds,300.0,19.0 %,Lagrange,temporal,300.0,10.9 %,1.0,$\blacksquare$,$\blacksquare$,{},$\blacksquare$,$\blacksquare$
Decomposing a renewable energy design and dispatch model,\cite{Wales2024},{REopt (CHP, TES)-2p},BU-CI-SHEX.te,{227,800 to 403,100 variables},AMPL,CPLEX 12,Optimality Gap,Runtime in seconds,300.0,9.7 %,Lagrange,temporal,300.0,3.2 %,1.0,$\blacksquare$,$\blacksquare$,{},$\blacksquare$,$\blacksquare$
Decomposing a renewable energy design and dispatch model,\cite{Wales2024},{REopt (CHP, TES)-5p-mT},BU-CI-SHEX.te,{227,800 to 403,100 variables},AMPL,CPLEX 12,Optimality Gap,Runtime in seconds,300.0,75.9 %,Lagrange,temporal,300.0,3.8 %,1.0,$\blacksquare$,$\blacksquare$,{},$\blacksquare$,$\blacksquare$
Decomposing a renewable energy design and dispatch model,\cite{Wales2024},{REopt (CHP, TES)-5p},BU-CI-SHEX.te,{227,800 to 403,100 variables},AMPL,CPLEX 12,Optimality Gap,Runtime in seconds,300.0,13.3 %,Lagrange,temporal,300.0,4.2 %,1.0,$\blacksquare$,$\blacksquare$,{},$\blacksquare$,$\blacksquare$
Decomposing a renewable energy design and dispatch model,\cite{Wales2024},{REopt (CHP, TES)-mT},BU-CI-SHEX.te,{227,800 to 403,100 variables},AMPL,CPLEX 12,Optimality Gap,Runtime in seconds,300.0,41.8 %,Lagrange,temporal,300.0,5.0 %,1.0,$\blacksquare$,$\blacksquare$,{},$\blacksquare$,$\blacksquare$
Decomposing a renewable energy design and dispatch model,\cite{Wales2024},REopt (CHP)-5p,BU-CI-SHEX.te,{227,800 to 403,100 variables},AMPL,CPLEX 12,Optimality Gap,Runtime in seconds,300.0,10.1 %,Lagrange,temporal,300.0,3.9 %,1.0,$\blacksquare$,$\blacksquare$,{},$\blacksquare$,$\blacksquare$
Decomposing a renewable energy design and dispatch model,\cite{Wales2024},{REopt (CHP, PV, BES)-5b-5p},BU-CI-SHEX.te,{227,800 to 403,100 variables},AMPL,CPLEX 12,Optimality Gap,Runtime in seconds,300.0,88.7 %,Lagrange,temporal,183.0,0.5 %,1.6,$\blacksquare$,$\blacksquare$,{},$\blacksquare$,$\blacksquare$
An Integrated Progressive Hedging and Benders Decomposition With Multiple Master Method to Solve the Brazilian Generation Expansion Problem,\cite{Soares2022},Custom GEP Model 20S,BU-CIS-EX.st,-,-,Xpress 34,Optimality Gap,Runtime in minutes,59.0,0.1 %,PH-Subgradient Multi-Master Benders,scenarios,46.0,0.1 %,1.3,$\blacksquare$,$\blacksquare$,$\blacksquare$,{},$\blacksquare$
Planning energy storage and photovoltaic panels for demand response with heating ventilation and air conditioning systems,\cite{Alhaider2018},Custom HVAC-BESS 600S,BU-CIS-SHEX.t,-,MATLAB,CPLEX -,Optimality Gap,Runtime in minutes,1440.0,-,Benders,variable types,35.0,0 %,41.1,$\blacksquare$,$\blacksquare$,$\blacksquare$,$\blacksquare$,$\blacksquare$
{Stabilized Benders decomposition for energy planning under climate uncertainty
},\cite{Gke2024},AnyMOD EuSysMod,BU-CS-EX.st,{4,577,298 variables},JuMP,Gurobi 10,Convergence,Runtime in seconds,1.6,yes,Benders,scenarios,7.6,yes,0.2,$\blacksquare$,{},$\blacksquare$,{},$\blacksquare$
{Stabilized Benders decomposition for energy planning under climate uncertainty
},\cite{Gke2024},AnyMOD EuSysMod,BU-CS-EX.st,{9,155,426 variables},JuMP,Gurobi 10,Convergence,Runtime in seconds,3.0,yes,Benders,scenarios,1.9,yes,1.5,$\blacksquare$,{},$\blacksquare$,{},$\blacksquare$
{Stabilized Benders decomposition for energy planning under climate uncertainty
},\cite{Gke2024},AnyMOD EuSysMod,BU-CS-EX.st,{13,732,768 variables},JuMP,Gurobi 10,Convergence,Runtime in seconds,4.7,yes,Benders,scenarios,6.8,yes,0.7,$\blacksquare$,{},$\blacksquare$,{},$\blacksquare$
{Stabilized Benders decomposition for energy planning under climate uncertainty
},\cite{Gke2024},AnyMOD EuSysMod,BU-CS-EX.st,{18,310,502 variables},JuMP,Gurobi 10,Convergence,Runtime in seconds,7.1,yes,Benders,scenarios,3.9,yes,1.8,$\blacksquare$,{},$\blacksquare$,{},$\blacksquare$
{Stabilized Benders decomposition for energy planning under climate uncertainty
},\cite{Gke2024},AnyMOD EuSysMod,BU-CS-EX.st,{22,887,920 variables},JuMP,Gurobi 10,Convergence,Runtime in seconds,10.3,yes,Benders,scenarios,6.6,yes,1.6,$\blacksquare$,{},$\blacksquare$,{},$\blacksquare$
{Stabilized Benders decomposition for energy planning under climate uncertainty
},\cite{Gke2024},AnyMOD EuSysMod,BU-CS-EX.st,{27,467,018 variables},JuMP,Gurobi 10,Convergence,Runtime in seconds,22.0,yes,Benders,scenarios,5.2,yes,4.3,$\blacksquare$,{},$\blacksquare$,{},$\blacksquare$
{Stabilized Benders decomposition for energy planning under climate uncertainty
},\cite{Gke2024},AnyMOD EuSysMod,BU-CS-EX.st,{32,044,314 variables},JuMP,Gurobi 10,Convergence,Runtime in seconds,20.2,yes,Benders,scenarios,5.9,yes,3.4,$\blacksquare$,{},$\blacksquare$,{},$\blacksquare$
{Stabilized Benders decomposition for energy planning under climate uncertainty
},\cite{Gke2024},AnyMOD EuSysMod,BU-CS-EX.st,{36,621,118 variables},JuMP,Gurobi 10,Convergence,Runtime in seconds,19.2,yes,Benders,scenarios,5.8,yes,3.3,$\blacksquare$,{},$\blacksquare$,{},$\blacksquare$
\end{filecontents*}
\csvreader[
            on column count error,
            head to column names,
            respect percent,
            respect underscore
        ]{results/DataExtraction-MethodBenchmark.csv}{}
        {\\ \hline \Citation & \ESOMName & \ESOMType & \ESOMSize & \AML & \BaselineSolver & \PerformanceMetric & \QualityMetric & \BaselinePerformance & \BaselineQuality & \DecompositionMethod & \DecompositionDimensions & \DecompositionPerformance & \DecompositionQuality}
        \\ \hline
    \end{longtable}
    \normalsize
    \end{landscape}

    \begin{landscape}
    \tiny
    \begin{longtable}{|
    m{.13\linewidth}|
    m{.08\linewidth}|
    m{.05\linewidth}|
    m{.05\linewidth}|
    m{.06\linewidth}|
    m{.06\linewidth}|
    m{.07\linewidth}|
    m{.06\linewidth}|
    m{.07\linewidth}|
    m{.04\linewidth}|
    m{.04\linewidth}|
    m{.05\linewidth}|}
        \caption{Scaling Benchmarks. \textit{Abbreviations: \acf{ESOM}, \acf{MEnv}, \acf{PMet}, \acf{QMet}, \acf{DecSys}, \acf{PVal}, \acf{QVal}}}
        \label{tab:results-scaling}
        \\ \hline
        \cellcolor{\tableHeaderColor}\textbf{Study} & 
        \cellcolor{\tableHeaderColor}\textbf{\acs{ESOM} Name} & 
        \cellcolor{\tableHeaderColor}\textbf{\acs{ESOM} Type} &
        \cellcolor{\tableHeaderColor}\textbf{\acs{ESOM} Size} &
        \cellcolor{\tableHeaderColor}\textbf{\acs{MEnv}} &
        \cellcolor{\tableHeaderColor}\textbf{\acs{PMet}} &
        \cellcolor{\tableHeaderColor}\textbf{\acs{QMet}} &
        \cellcolor{\tableHeaderColor}\textbf{\acs{DecSys} Method} &
        \cellcolor{\tableHeaderColor}\textbf{Decomposed Dimensions} &
        \cellcolor{\tableHeaderColor}\textbf{\acs{DecSys} \acs{PVal}} &
        \cellcolor{\tableHeaderColor}\textbf{\acs{DecSys} \acs{QVal}} &
        \cellcolor{\tableHeaderColor}\textbf{Cores}
        \begin{filecontents*}{results/DataExtraction-ScalingBenchmark.csv}
PublicationTitle,Citation,ESOMName,ESOMType,ESOMSize,AML,QualityMetric,PerformanceMetric,DecompositionMethod,DecompositionDimensions,DecompositionPerformance,DecompositionQuality,ParallelProcesses,Continuity,Integrality,Stochasticity,SHDecisions,EXDecisions
Short-term Hydrothermal Generation Scheduling Using a Parallelized Stochastic Mixed-integer Linear Programming Algorithm,\cite{Gil2016},Custom 12 Scenario STHTGS,BU-CIS-SH.ste,{5,124,120 variables},FORTRAN,Convergence,Runtime in hours,Progressive Hedging,scenarios,39.0,yes,2,$\blacksquare$,$\blacksquare$,$\blacksquare$,$\blacksquare$,
Short-term Hydrothermal Generation Scheduling Using a Parallelized Stochastic Mixed-integer Linear Programming Algorithm,\cite{Gil2016},Custom 12 Scenario STHTGS,BU-CIS-SH.ste,{5,124,120 variables},FORTRAN,Convergence,Runtime in hours,Progressive Hedging,scenarios,23.0,yes,4,$\blacksquare$,$\blacksquare$,$\blacksquare$,$\blacksquare$,
Short-term Hydrothermal Generation Scheduling Using a Parallelized Stochastic Mixed-integer Linear Programming Algorithm,\cite{Gil2016},Custom 12 Scenario STHTGS,BU-CIS-SH.ste,{5,124,120 variables},FORTRAN,Convergence,Runtime in hours,Progressive Hedging,scenarios,17.0,yes,8,$\blacksquare$,$\blacksquare$,$\blacksquare$,$\blacksquare$,
Short-term Hydrothermal Generation Scheduling Using a Parallelized Stochastic Mixed-integer Linear Programming Algorithm,\cite{Gil2016},Custom 12 Scenario STHTGS,BU-CIS-SH.ste,{5,124,120 variables},FORTRAN,Convergence,Runtime in hours,Progressive Hedging,scenarios,9.0,yes,12,$\blacksquare$,$\blacksquare$,$\blacksquare$,$\blacksquare$,
Short-term Hydrothermal Generation Scheduling Using a Parallelized Stochastic Mixed-integer Linear Programming Algorithm,\cite{Gil2016},Custom 12 Scenario STHTGS,BU-CIS-SH.ste,{5,124,120 variables},FORTRAN,Convergence,Runtime in hours,Distributed Progressive Hedging,scenarios,5.0,yes,24,$\blacksquare$,$\blacksquare$,$\blacksquare$,$\blacksquare$,
High-performance computing based fully parallel security-constrained unit commitment with dispatchable transmission network,\cite{Gong2019},Modified IEEE 118-Bus 60sw-30c,BU-CI-SH.st,{386,496 variables},-,Convergence,Runtime in hours,Lagrange APP-AL,{operational, temporal, spatial, technological},15.7,yes,1,$\blacksquare$,$\blacksquare$,{},$\blacksquare$,
High-performance computing based fully parallel security-constrained unit commitment with dispatchable transmission network,\cite{Gong2019},Modified IEEE 118-Bus 60sw-30c,BU-CI-SH.st,{386,496 variables},-,Convergence,Runtime in hours,Lagrange APP-AL,{operational, temporal, spatial, technological},9.3,yes,2,$\blacksquare$,$\blacksquare$,{},$\blacksquare$,
High-performance computing based fully parallel security-constrained unit commitment with dispatchable transmission network,\cite{Gong2019},Modified IEEE 118-Bus 60sw-30c,BU-CI-SH.st,{386,496 variables},-,Convergence,Runtime in hours,Lagrange APP-AL,{operational, temporal, spatial, technological},4.9,yes,4,$\blacksquare$,$\blacksquare$,{},$\blacksquare$,
High-performance computing based fully parallel security-constrained unit commitment with dispatchable transmission network,\cite{Gong2019},Modified IEEE 118-Bus 60sw-30c,BU-CI-SH.st,{386,496 variables},-,Convergence,Runtime in hours,Lagrange APP-AL,{operational, temporal, spatial, technological},2.6,yes,8,$\blacksquare$,$\blacksquare$,{},$\blacksquare$,
High-performance computing based fully parallel security-constrained unit commitment with dispatchable transmission network,\cite{Gong2019},Modified IEEE 118-Bus 60sw-30c,BU-CI-SH.st,{386,496 variables},-,Convergence,Runtime in hours,Lagrange APP-AL,{operational, temporal, spatial, technological},1.3,yes,16,$\blacksquare$,$\blacksquare$,{},$\blacksquare$,
High-performance computing based fully parallel security-constrained unit commitment with dispatchable transmission network,\cite{Gong2019},Modified IEEE 118-Bus 60sw-30c,BU-CI-SH.st,{386,496 variables},-,Convergence,Runtime in hours,Lagrange APP-AL,{operational, temporal, spatial, technological},0.7,yes,32,$\blacksquare$,$\blacksquare$,{},$\blacksquare$,
High-performance computing based fully parallel security-constrained unit commitment with dispatchable transmission network,\cite{Gong2019},Modified IEEE 118-Bus 60sw-30c,BU-CI-SH.st,{386,496 variables},-,Convergence,Runtime in hours,Lagrange APP-AL,{operational, temporal, spatial, technological},0.4,yes,64,$\blacksquare$,$\blacksquare$,{},$\blacksquare$,
High-performance computing based fully parallel security-constrained unit commitment with dispatchable transmission network,\cite{Gong2019},Modified IEEE 118-Bus 60sw-30c,BU-CI-SH.st,{386,496 variables},-,Convergence,Runtime in hours,Lagrange APP-AL,{operational, temporal, spatial, technological},0.3,yes,128,$\blacksquare$,$\blacksquare$,{},$\blacksquare$,
High-performance computing based fully parallel security-constrained unit commitment with dispatchable transmission network,\cite{Gong2019},Modified IEEE 118-Bus 60sw-30c,BU-CI-SH.st,{386,496 variables},-,Convergence,Runtime in hours,Lagrange APP-AL,{operational, temporal, spatial, technological},0.3,yes,159,$\blacksquare$,$\blacksquare$,{},$\blacksquare$,
High-performance computing based fully parallel security-constrained unit commitment with dispatchable transmission network,\cite{Gong2019},Modified IEEE 118-Bus 60sw-30c,BU-CI-SH.st,{386,496 variables},-,Convergence,Runtime in hours,Lagrange APP-AL,{operational, temporal, spatial, technological},0.3,yes,179,$\blacksquare$,$\blacksquare$,{},$\blacksquare$,
High-performance computing based fully parallel security-constrained unit commitment with dispatchable transmission network,\cite{Gong2019},Modified IEEE 118-Bus 60sw-30c,BU-CI-SH.st,{386,496 variables},-,Convergence,Runtime in hours,Lagrange APP-AL,{operational, temporal, spatial, technological},0.2,yes,199,$\blacksquare$,$\blacksquare$,{},$\blacksquare$,
High-performance computing based fully parallel security-constrained unit commitment with dispatchable transmission network,\cite{Gong2019},IEEE 1168-Bus 60sw-30c,BU-CI-SH.st,{734,256 variables},-,Convergence,Runtime in hours,Lagrange APP-AL,{operational, temporal, spatial, technological},22.2,yes,1,$\blacksquare$,$\blacksquare$,{},$\blacksquare$,
High-performance computing based fully parallel security-constrained unit commitment with dispatchable transmission network,\cite{Gong2019},IEEE 1168-Bus 60sw-30c,BU-CI-SH.st,{734,256 variables},-,Convergence,Runtime in hours,Lagrange APP-AL,{operational, temporal, spatial, technological},11.6,yes,2,$\blacksquare$,$\blacksquare$,{},$\blacksquare$,
High-performance computing based fully parallel security-constrained unit commitment with dispatchable transmission network,\cite{Gong2019},IEEE 1168-Bus 60sw-30c,BU-CI-SH.st,{734,256 variables},-,Convergence,Runtime in hours,Lagrange APP-AL,{operational, temporal, spatial, technological},6.3,yes,4,$\blacksquare$,$\blacksquare$,{},$\blacksquare$,
High-performance computing based fully parallel security-constrained unit commitment with dispatchable transmission network,\cite{Gong2019},IEEE 1168-Bus 60sw-30c,BU-CI-SH.st,{734,256 variables},-,Convergence,Runtime in hours,Lagrange APP-AL,{operational, temporal, spatial, technological},3.8,yes,8,$\blacksquare$,$\blacksquare$,{},$\blacksquare$,
High-performance computing based fully parallel security-constrained unit commitment with dispatchable transmission network,\cite{Gong2019},IEEE 1168-Bus 60sw-30c,BU-CI-SH.st,{734,256 variables},-,Convergence,Runtime in hours,Lagrange APP-AL,{operational, temporal, spatial, technological},2.0,yes,16,$\blacksquare$,$\blacksquare$,{},$\blacksquare$,
High-performance computing based fully parallel security-constrained unit commitment with dispatchable transmission network,\cite{Gong2019},IEEE 1168-Bus 60sw-30c,BU-CI-SH.st,{734,256 variables},-,Convergence,Runtime in hours,Lagrange APP-AL,{operational, temporal, spatial, technological},1.2,yes,32,$\blacksquare$,$\blacksquare$,{},$\blacksquare$,
High-performance computing based fully parallel security-constrained unit commitment with dispatchable transmission network,\cite{Gong2019},IEEE 1168-Bus 60sw-30c,BU-CI-SH.st,{734,256 variables},-,Convergence,Runtime in hours,Lagrange APP-AL,{operational, temporal, spatial, technological},1.2,yes,64,$\blacksquare$,$\blacksquare$,{},$\blacksquare$,
High-performance computing based fully parallel security-constrained unit commitment with dispatchable transmission network,\cite{Gong2019},IEEE 1168-Bus 60sw-30c,BU-CI-SH.st,{734,256 variables},-,Convergence,Runtime in hours,Lagrange APP-AL,{operational, temporal, spatial, technological},1.0,yes,128,$\blacksquare$,$\blacksquare$,{},$\blacksquare$,
High-performance computing based fully parallel security-constrained unit commitment with dispatchable transmission network,\cite{Gong2019},IEEE 1168-Bus 60sw-30c,BU-CI-SH.st,{734,256 variables},-,Convergence,Runtime in hours,Lagrange APP-AL,{operational, temporal, spatial, technological},0.8,yes,159,$\blacksquare$,$\blacksquare$,{},$\blacksquare$,
High-performance computing based fully parallel security-constrained unit commitment with dispatchable transmission network,\cite{Gong2019},IEEE 1168-Bus 60sw-30c,BU-CI-SH.st,{734,256 variables},-,Convergence,Runtime in hours,Lagrange APP-AL,{operational, temporal, spatial, technological},0.8,yes,179,$\blacksquare$,$\blacksquare$,{},$\blacksquare$,
High-performance computing based fully parallel security-constrained unit commitment with dispatchable transmission network,\cite{Gong2019},IEEE 1168-Bus 60sw-30c,BU-CI-SH.st,{734,256 variables},-,Convergence,Runtime in hours,Lagrange APP-AL,{operational, temporal, spatial, technological},0.8,yes,199,$\blacksquare$,$\blacksquare$,{},$\blacksquare$,
Distributed Price-Directive Decomposition Applications in Power Systems Operations,\cite{Sundarraj1995},Custom Power Network 17G,BU-C-SH.st,-,FORTRAN,Convergence,Runtime in seconds,Dantzig-Wolfe,topological,5.5,yes,2,$\blacksquare$,{},,$\blacksquare$,
Distributed Price-Directive Decomposition Applications in Power Systems Operations,\cite{Sundarraj1995},Custom Power Network 17G,BU-C-SH.st,-,FORTRAN,Convergence,Runtime in seconds,Dantzig-Wolfe,topological,4.3,yes,5,$\blacksquare$,{},,$\blacksquare$,
Distributed Price-Directive Decomposition Applications in Power Systems Operations,\cite{Sundarraj1995},Custom Power Network 17G,BU-C-SH.st,-,FORTRAN,Convergence,Runtime in seconds,Dantzig-Wolfe,topological,4.5,yes,10,$\blacksquare$,{},,$\blacksquare$,
Distributed Price-Directive Decomposition Applications in Power Systems Operations,\cite{Sundarraj1995},Custom Power Network 17G,BU-C-SH.st,-,FORTRAN,Convergence,Runtime in seconds,Dantzig-Wolfe,topological,5.5,yes,15,$\blacksquare$,{},,$\blacksquare$,
Distributed Price-Directive Decomposition Applications in Power Systems Operations,\cite{Sundarraj1995},Custom Power Network 68G,BU-C-SH.st,-,FORTRAN,Convergence,Runtime in seconds,Dantzig-Wolfe,topological,19.8,yes,5,$\blacksquare$,{},,$\blacksquare$,
Distributed Price-Directive Decomposition Applications in Power Systems Operations,\cite{Sundarraj1995},Custom Power Network 68G,BU-C-SH.st,-,FORTRAN,Convergence,Runtime in seconds,Dantzig-Wolfe,topological,10.2,yes,10,$\blacksquare$,{},,$\blacksquare$,
Distributed Price-Directive Decomposition Applications in Power Systems Operations,\cite{Sundarraj1995},Custom Power Network 68G,BU-C-SH.st,-,FORTRAN,Convergence,Runtime in seconds,Dantzig-Wolfe,topological,8.0,yes,15,$\blacksquare$,{},,$\blacksquare$,
Distributed Price-Directive Decomposition Applications in Power Systems Operations,\cite{Sundarraj1995},Custom Power Network 68G,BU-C-SH.st,-,FORTRAN,Convergence,Runtime in seconds,Dantzig-Wolfe,topological,8.6,yes,20,$\blacksquare$,{},,$\blacksquare$,
Distributed Price-Directive Decomposition Applications in Power Systems Operations,\cite{Sundarraj1995},Custom Power Network 68G,BU-C-SH.st,-,FORTRAN,Convergence,Runtime in seconds,Dantzig-Wolfe,topological,8.9,yes,25,$\blacksquare$,{},,$\blacksquare$,
Distributed Price-Directive Decomposition Applications in Power Systems Operations,\cite{Sundarraj1995},Custom Power Network 68G,BU-C-SH.st,-,FORTRAN,Convergence,Runtime in seconds,Dantzig-Wolfe,topological,10.1,yes,30,$\blacksquare$,{},,$\blacksquare$,
Distributed Price-Directive Decomposition Applications in Power Systems Operations,\cite{Sundarraj1995},Custom Power Network 119G,BU-C-SH.st,-,FORTRAN,Convergence,Runtime in seconds,Dantzig-Wolfe,topological,46.1,yes,5,$\blacksquare$,{},,$\blacksquare$,
Distributed Price-Directive Decomposition Applications in Power Systems Operations,\cite{Sundarraj1995},Custom Power Network 119G,BU-C-SH.st,-,FORTRAN,Convergence,Runtime in seconds,Dantzig-Wolfe,topological,15.9,yes,10,$\blacksquare$,{},,$\blacksquare$,
Distributed Price-Directive Decomposition Applications in Power Systems Operations,\cite{Sundarraj1995},Custom Power Network 119G,BU-C-SH.st,-,FORTRAN,Convergence,Runtime in seconds,Dantzig-Wolfe,topological,12.8,yes,15,$\blacksquare$,{},,$\blacksquare$,
Distributed Price-Directive Decomposition Applications in Power Systems Operations,\cite{Sundarraj1995},Custom Power Network 119G,BU-C-SH.st,-,FORTRAN,Convergence,Runtime in seconds,Dantzig-Wolfe,topological,11.9,yes,20,$\blacksquare$,{},,$\blacksquare$,
Distributed Price-Directive Decomposition Applications in Power Systems Operations,\cite{Sundarraj1995},Custom Power Network 119G,BU-C-SH.st,-,FORTRAN,Convergence,Runtime in seconds,Dantzig-Wolfe,topological,11.9,yes,25,$\blacksquare$,{},,$\blacksquare$,
Distributed Price-Directive Decomposition Applications in Power Systems Operations,\cite{Sundarraj1995},Custom Power Network 119G,BU-C-SH.st,-,FORTRAN,Convergence,Runtime in seconds,Dantzig-Wolfe,topological,12.6,yes,30,$\blacksquare$,{},,$\blacksquare$,
\end{filecontents*}
\csvreader[
            on column count error,
            head to column names,
            respect percent,
            respect underscore            
        ]{results/DataExtraction-ScalingBenchmark.csv}{}
        {\\ \hline \Citation & \ESOMName & \ESOMType & \ESOMSize & \AML & \PerformanceMetric &  \QualityMetric & \DecompositionMethod & \DecompositionDimensions & \DecompositionPerformance & \DecompositionQuality & \ParallelProcesses}
        \\ \hline        
    \end{longtable}
    \normalsize
    \end{landscape}

\end{document}

%% file: sections/01-introduction.tex
\section{Introduction}
\label{sec:introduction}

    Energy supply systems are currently undergoing structural and regulatory changes while scaling up. These factors are reflected by increasing dimensionality and connectivity of mathematical energy system optimization models. Numerical methods for solving optimization models have significantly improved over the last decades \cite{koch2022progress}, likewise has the performance of modern computation processors \cite{gonzalez2019trends}.

    Performance gains had primarily been derived from enhancements in the efficiency of sequential processing, both in algorithms and hardware. In order to advance beyond the current barriers of sequential processing, computing hardware has developed towards  parallel processing \cite{millett2011future}. Iterative optimization algorithms may catch up by means of a similar approach \cite{zhou2023review}. Any further performance improvements may be derived from novel computation processors \cite{shalf2020future} or algorithm engineering \cite{sanders2009algorithm}.

    Efforts to speed up solving linear programs with and without integrality constraints have led to the exploration of parallelization strategies for optimization algorithms as well as applications of decomposition techniques for optimization models. Parallelization can be done on the functional level by decomposing the algorithm into independent tasks or parallelized units. Another approach is the decomposition of the problem within its domain. Previous work that surveys these attempts for global optimization and methods tailored to energy system optimization is lined out in the following.

    \subsubsection*{Parallelized Exact Optimization}
        We are going to start on the simplex algorithm \cite{dantzig1955generalized}, an iterative method traversing the edges of the feasible region towards the optimum with numerous, computationally inexpensive steps. Development of parallelized versions will not pay off in most of the cases as a review on its parallelization shows \cite{hall2010towards}. Parallelized simplex methods could not outperform the corresponding serial implementations.
    
        Next, we consider interior point methods \cite{dikin1967iterative}, which converge to the optimum while traversing the interior of the feasible region with fewer, computationally expensive steps. Given a required substructure of the problem's formulation, interior point methods benefit from the exploitation of this substructure on each iteration \cite{gondzio2003parallel,gondzio2005direct,gondzio2012interior}. Algorithms of this class also extend to nonlinear problems.

        We proceed on the Branch-And-Bound method \cite{land1960automatic}, which is a systematic tree search strategy that renders algorithms which incorporate integrality conditions. Parallelization can be effective for Branch-And-Bound based algorithms as shown in the last comprehensive survey on this topic \cite{gendron1994parallel}. As the number of available compute nodes has significantly grown, the limits of parallelization in the context of Branch-And-Bound algorithms has been assessed \cite{koch2012could}, concluding that the solution of node relaxations plays a major role. Parallelization strategies need to deal with additional challenges such as selection rules and load balancing to be efficient as pointed out by \citet{herrera2017parallel}, who investigate how the implementation framework influences the algorithm's performance.

        Finally, fist-order methods \cite{cauchy1847methode} are considered. They only use derivatives not higher than of first order to iteratively move along the gradient towards the optimum while modifying their update steps or objective function to incorporate constraints \cite{beck2017first}. A review by \citet{liu2022survey} highlights the use of distributed environments for gradient-based methods focusing on nonlinear optimization \cite{liu2022survey}. The convergence rate could be improved for a hybrid gradient method that alternates between the primal and dual formulation \cite{Zhu2008AnEP}, the performance of which has been significantly increased by \citet{Applegate2021} and parallelized by \citet{applegate2025pdlp}.
        
        Beyond the parallelization of solving algorithms, optimization models may also be decomposed. Given a certain substructure, model decomposition in the context of linear optimization leads to different hierarchical decomposition methods. Exact model decomposition techniques can keep high accuracy and naturally distribute on modern high performance computing environments \cite{karbowski2015decomposition}. 

        Given the previous work, we can conclude that simplex methods have limited potential for concurrent computation. Interior point methods mostly benefit from data parallelism given a substructure in the model. Branching methods offer better opportunities for task parallelism while being challenging for computational load balancing.
 
    \subsubsection*{Energy Systems Decomposition And Parallelization}
        The existing literature for energy system optimization reviews either decomposition or parallelization methods independently. A research article by \citet{Sagastizbal2012DivideTC} explores various decomposition techniques to address the growing complexity of energy systems. It investigates a set of decomposition methods on six prototypical examples, providing qualitative assessments of the selected methods and including two case studies. Another survey targeting power systems \cite{Molzahn2017ASO} reviews techniques to implement distributed optimization algorithms for either linear or convex-nonlinear or nonconvex optimal power flow models. The authors categorize the methods into either augmented Lagrangian decomposition or decentralized solution of the Karush–Kuhn–Tucker optimality conditions. A systematic evaluation by \citet{Cao2019ClassificationAE} reviews several aggregation and decomposition methods for models based on the REMix-Framework. Their evaluation covers aggregation methods as well as two heuristic approaches which temporally decompose the system at reduced resolution. Furthermore, \citet{Rodriguez2021ARO} cover parallel heterogeneous computing techniques for optimization and analysis of power systems. Another review by \citet{al2022high}, focused on electrical energy system optimization, gives an overview on the different types of hardware that allow for parallelization of the solution procedures.      

\subsubsection*{Scope of the review}

    The previous work has focused either on decomposition methods to manage model complexity or on parallelization to improve computational efficiency without the integration of both. This study addresses that gap by providing the first comprehensive and traceable survey of parallelized decomposition approaches benchmarked within the context of linear energy system models. In this context, we classify the associated benchmark models and examine software systems that are particularly well-suited for supporting such parallel approaches. This allows for a structured comparison of parallelized methods suitable for most of the relevant energy system models which are linear.
    
    Breaking down large optimization problems into smaller, independent sub-problems, decomposition methods promise performance improvements as they terminate earlier or run in parallel. Therefore, we are stating the following questions:
    
    \begin{flushleft}
    \begin{enumerate}
        \item Which classes of \acp{ESOM} can be defined?
        \item Which classes of parallelized decomposition methods are frequently employed?
        \item How do the parallelized decomposition methods perform on the different model classes?
    \end{enumerate}
    \end{flushleft}

    The following \autoref{sec:theory} yields an overview of the basic theory and terminology. The subsequent \autoref{sec:methods} introduces the methods we employed to review the literature on parallelized decomposition in energy system optimization. \autoref{sec:esoms} covers \acp{ESOM} and their properties. In \autoref{sec:parallelization}, decomposition methods are introduced as a means to parallelization. \autoref{sec:comparison} yields a comparison of these methods with respect to our stated questions and formulates recommendations for conducting benchmark studies. In the last \autoref{sec:conclusion} we conclude on the results about decomposition methods for various energy system models as a way to improve computational performance.

    The mathematical notation in this publication follows part two of the ISO 80000 standard: Matrices are written with bold italic capital letters and their elements with thin italic lowercase letters. Vectors are written as bold italic lowercase letters and scalars are thin italic lowercase letters. All acronyms are listed in \autoref{sec:acronyms}.

%% file: sections/02-theory.tex
\section{Theory}
\label{sec:theory}

In this section, we are first going to give a primer on polytopes as to introduce some basic terminology, for reference compare further \cite{villavicencio2024polytopes} and \cite{ziegler2012lectures}.

A polytope is a bounded polyhedron and a polyhedron is the intersection of finitely many closed halfspaces and therefore is always a convex set. Every polyhedron has two equal representations, either as the intersection of its determining halfspaces, referred to as $\mathcal{H}$-representation, or as the Minkowski sum of the convex hull of its vertices and the conical hull of its rays, referred to as $\mathcal{V}$-representation. The equality of those representations is stated by the Weyl-Minkowski theorem \cite{weyl1934elementare}. The convex hull of a set of points is the set of all convex combinations of those points. The conical hull of a set of points is the set of all affine combinations of those points. A linear combination, i.e. the weighted sum, of a set of points is conical if and only if all coefficients are non-negative. If all coefficients add up to one, it is affine. A convex combination is defined as a linear combination that is both conical and affine. Furthermore, the product of two polytopes is the Cartesian product of their defining sets and results in another polytope. Given an irredundant $\mathcal{H}$-representation of a $d$-dimensional convex polytope, a $k$-face is the set of points which fulfill $d-k$ of the determining inequalities as an equality.

Next, we yield some terminology for parallel computing, for reference compare further \cite{padua2011encyclopedia} and \cite{lin2008principles}.

A process is a program being executed with its assigned system resources and its context. A parallel program simultaneously performs multiple processes. Given a work load that has been decomposed and assigned to several processes, a system that allows those processes to share the same primary memory is a shared-memory parallel system, whereas a system in which the processes exchange information only via explicit communication is a distributed-memory system. If the input data can be partitioned in a highly granular way such that the same operations are executed in parallel on the different partitions, we call this data-parallelism. If different blocks of operations, the tasks, are executed on the same or on different partitions of the input, we refer to it as task-parallelism. In a parallel system multiple processes might request access to a shared resource which leads to contention. A multicore system also needs to keep its memory state coherent which introduces additional delay. A widely used classification scheme for parallel computer architectures is Flynn’s taxonomy \cite{flynn1972some}, classifying by microprocessor-level instruction stream and data stream processing, defining the following catagories: SISD (Single Instruction, Single Data), operating one instruction on a single data stream, possibly taking advantage of instruction-level parallelism within the instruction stream, e.g. pipelining. SIMD (Single Instruction, Multiple Data), applying an instruction on multiple data streams in parallel, e.g. array processors.  MISD (Multiple Instruction, Single Data), processing one data stream on different processing units. MIMD (Multiple Instruction, Multiple Data), performing different instructions on multiple data streams, e.g. multi-threaded and multi-core processors.

Lastly, a brief overview on computational performance analysis is given, for reference compare further \cite{liu2011software} and \cite{lilja2005measuring}.

A performance metric is a time, count or size value that measures the system's performance we are interested in, and should be linear, reliable, repeatable and consistent. A performance metric normalized to a time unit is referenced to as throughput. A benchmark system's speedup $s$ compared to a reference system is the ratio of its throughput $R$ and the reference system's throughput $R_{ref}$, i.e. $s=\frac{R}{R_{ref}}=\frac{T_{ref}}{T}$ with the benchmark system's runtime $T$ and the reference system's runtime $T_{ref}$.
Amdahl's law which has been derived from Amdahl's arguement \cite{amdahl1967validity} is given as $s=(f-\frac{1-f}{P})^{-1}$ for $P$ processors, with $f$ as the fraction that amounts to the not parallelizable part of the program. This relation assumes a fixed problem size and a variable number of parallel processors and is referred to as strong scaling. If both problem size and number of processors are variable, weak scaling is measured with a constant workload per processor, described by Gustafson's law \cite{gustafson1988reevaluating} as $s=f+P(1-f)$. A more detailed relation for the throughput $R(P)$ with $P$ processors considering the parallel system's contention level $\alpha$ and coherency delay $\beta$ is given by Gunther's law \cite{gunther1993simple} as $R(P)/R(1)=P\cdot(1+\alpha(P-1)+\beta P(P-1))^{-1}$. As complex computing system's are subject to performance variability, measurements are supposed to be sampled and given as a mean and its corresponding variability metric.

%% file: sections/03-methods.tex
\section{Review Methods}
\label{sec:methods}

    \begin{wrapfigure}{r}{0.5\textwidth}
        \includegraphics[width=\linewidth]{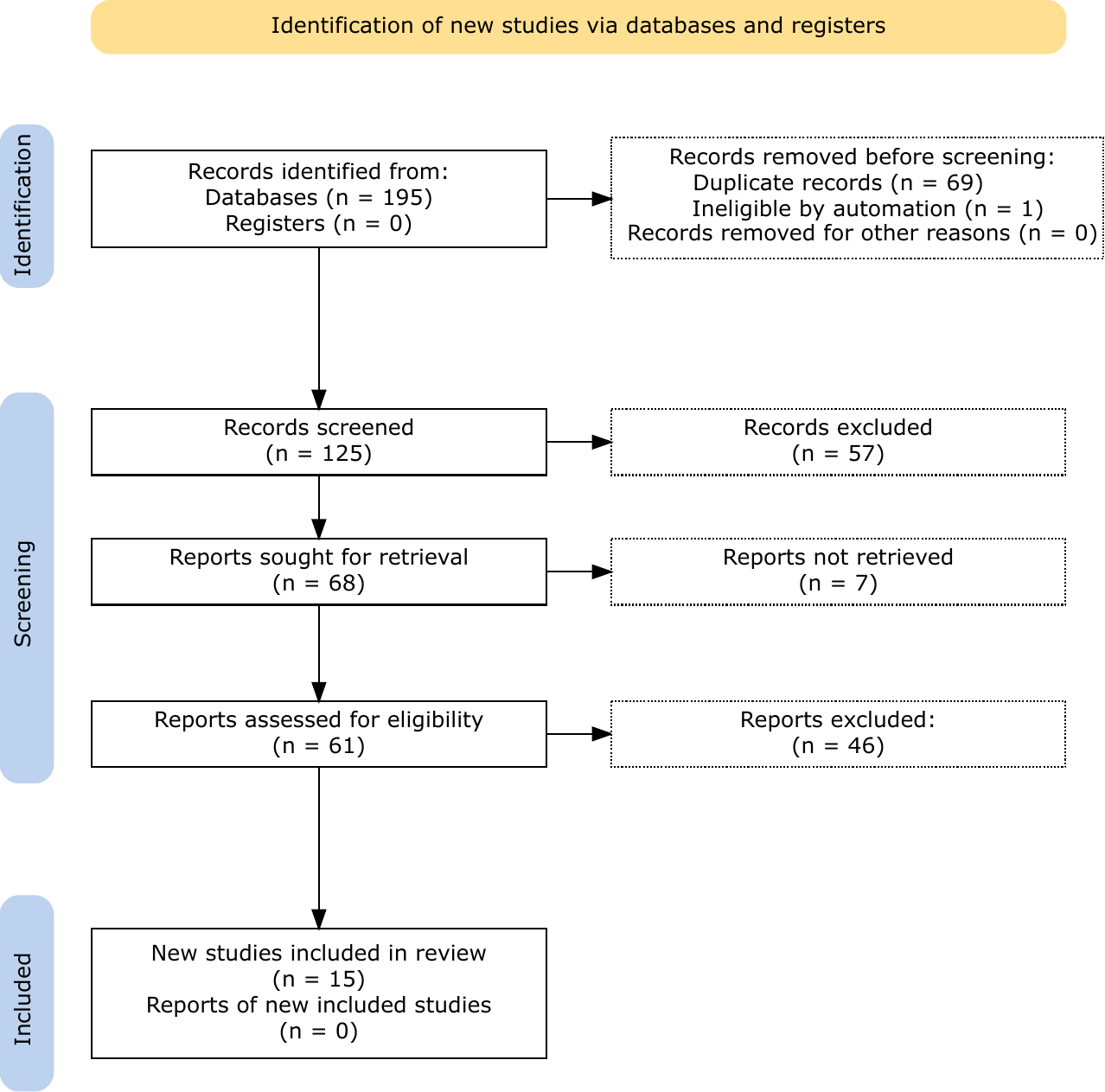}
        \caption{Number of records identified, included and excluded in the present review.}
        \label{fig:prisma-flow}
    \end{wrapfigure}

We conduct a systematic review on parallelized model decomposition strategies in the context of linear energy system optimization. For this, we collect, analyze and extract findings from the literature and sum up the interpretations. In order to make the work reproducible, the PRISMA statement \cite{Pagen71} is employed for tracing the review process. After retrieval, the records have been deduplicated by the Systematic Review Accelerator as it provides a traceable automatic procedure \cite{Forbes2024} as outlined in a publication on deduplication tools by \citet{guimaraes2022deduplicating}. The screening process has been done with Rayyan \cite{ouzzani2016rayyan}, which detected further duplicates. All records' abstracts have been screened for their relevance, i.e. a study that employs parallelization and uses a known decomposition method outlined in \autoref{sec:parallelization}. Publications on nonlinear models or (meta-)heuristics are excluded as well as methods such as multi-objective optimization, bi-level programming or equilibrium programs. 

While sections \ref{sec:esoms} and \ref{sec:parallelization} provide overviews of energy system optimisation models and parallel decomposition, the systematic review process described in this section was used for the results in \autoref{sec:comparison}. 

\subsection{Reporting Guideline}    
    Primarily targeted at meta-analyses and systematic reviews for evaluating health interventions, various extensions to the main PRISMA statement provide guidance for different types of systematic reviews. The guidelines help to clearly communicate how a systematic review was conducted, which methods were used, and which findings have been obtained. The PRISMA-S guideline \cite{rethlefsen2021prisma} includes 16 reporting items we are using to document the search strategy. This guideline is broadly applicable and therefore suitable for systematic reviews in a variety of fields.

\subsection{Preferred Reporting Items}
    The preferred reporting items of the selected guideline can be categorized into Information Sources, Search Strategies, Peer Review and Record Management. The detailed documentation of all items is shown in \autoref{tab:prisma-checklist}. The flow diagram in \autoref{fig:prisma-flow} shows the information flow through the different phases of the review. It provides an overview on the selection process, tracing the decisions made at each stage of the process. We queried the two literature databases Scopus and Web of Science to search for relevant publications. In both cases, we have searched for title, abstract, and author keywords. The Web of Science platform additionally includes terms generated from the titles of referenced papers, which are processed using a ranking algorithm \cite{Garfield1993KeyWordsPD}. We are seeking optimization models which are only linear and large in scale, excluding all nonlinear models. Among the optimization models, we narrow down the topic to \acp{ESOM}. This accounts for the majority of energy systems related publications. As real-time control models showed up frequently, they have been excluded. We want to retrieve only publications that focus on improving computational performance or tractability of the models. This term ensures that all studies are related to any kind of performance improvement that is typically found in the context of high performance computing or parallelized computing is captured, possibly sequential improvements, too, as we want to make sure that no study is lost if parallelization has been employed but not highlighted as the study's focus. Finally, the publications are supposed to speed up the solution process by any kind of decomposition. All model configurations were then examined in detail within the identified publications. This means that the number of configurations examined far exceeds the number of publications.
    
    \begin{table*}
    \caption{PRISMA-S 16-item checklist}
    \label{tab:prisma-checklist}
    \tiny
        \begin{tabular}{R{0.15\linewidth}R{0.25\linewidth}R{0.53\linewidth}}
        \hline
            \cellcolor{\tableHeaderColor}\textbf{Section} & \cellcolor{\tableHeaderColor}\textbf{Item} & \cellcolor{\tableHeaderColor}\textbf{Report} \\[1em] \hline
            \multicolumn{3}{l}{\textbf{INFORMATION SOURCES AND METHODS}} \\[0.5em] \hline
            Database name & The individual databases searched. & SCIE-EXPANDED, CPCI-S, BKCI-S, ESCI, SCOPUS \\ \hline
            Multi-database searching & Name of the platform searching databases  simultaneously. & Web Of Science: SCIE-EXPANDED, CPCI, BKCI-S, ESCI \\ \hline
            Study registries & List of the study registries searched. & none \\ \hline
            Online resources & Web search engines, web sites or other resource searched. & none \\ \hline
            Online resources and browsing & Online or print source purposefully searched or browsed. & none \\ \hline
            Citation searching & Cited references or citing references examined. & ~ \\ \hline
            Contacts & Publications by contacting authors or other experts. & ~ \\ \hline
            Other methods & Additional sources or search methods used. & none \\ \hline
            \multicolumn{3}{l}{\textbf{SEARCH STRATEGIES}} \\[0.5em] \hline
            Full search strategies  & Search strategies for each database and information source, exactly as run.  & 
            \vspace{-\topsep}
            \begin{lstlisting}[basicstyle=\tiny, xleftmargin=\dimexpr-\leftmargin-\leftmargini+1em, breaklines, breakatwhitespace, title=\raggedright{\textit{\tiny Scopus:}}, ]
                TITLE-ABS-KEY(large AND linear AND optim* AND (distributed-computing OR parallel OR hpc OR super-computing OR supercomputing OR cluster-computing OR clustercomputing OR solvable OR tractable OR feasible OR speed-up OR speedup OR convergence OR computation-time OR solution-time OR outperform OR performance) AND ((energy OR heat OR gas OR power OR electricity) W/0 (system OR network OR grid OR market)) AND (decompos*) AND NOT non-linear AND NOT nonlinear AND NOT control)
            \end{lstlisting}
            \vspace{-\topsep}
            \begin{lstlisting}[basicstyle=\tiny, xleftmargin=\dimexpr-\leftmargin-\leftmargini+1em, breaklines, breakatwhitespace, title=\raggedright{\textit{\tiny Web Of Science:}}]
                TS=(large AND linear NOT non-linear NOT nonlinear NOT control AND optim* AND (distributed-computing OR parallel OR hpc OR super-computing OR supercomputing OR cluster-computing OR clustercomputing OR solvable OR tractable OR feasible OR speed-up OR speedup OR convergence OR computation-time OR solution-time OR outperform OR performance) AND ((energy OR heat OR gas OR power OR electricity) NEAR/0 (system OR network OR grid OR market)) AND (decompos*))
            \end{lstlisting}
            \\ \hline
            Limits and restrictions & Limits or restrictions applied to a search. & none \\ \hline
            Search filters & Published search filters used (original or modified). & none \\ \hline
            Prior work & Search strategies from other literature reviews. & none \\ \hline
            Updates & Methods used to update the search. & none \\ \hline
            Dates of searches & For each search strategy, date of last search occurred. & 
            \begin{itemize}[leftmargin=*,noitemsep,topsep=0pt]
                \item Web Of Science: Jan. 14, 2025
                \item Scopus: Jan. 14, 2025
            \end{itemize}
            \\ \hline
            \multicolumn{3}{l}{\textbf{PEER REVIEW}} \\[0.5em] \hline
            Peer review & Description of any search peer review process.  & ~ \\ \hline
            \multicolumn{3}{l}{\textbf{MANAGING RECORDS}} \\[0.5em] \hline
            Total Records & Total number of records identified from all sources. & 195 \\ \hline
            Deduplication & Processes and software used to deduplicate records. & Export of full record RIS from each platform. Import to SRA Deduplicator, Focused-Mode Algorithm and export to RIS. Import of RIS into Rayyan, manual deduplication of remaining detections there.\\ \hline
        \end{tabular}
    \end{table*}

%% file: sections/04-energySystemModels.tex
\section{Energy System Optimization Models}
\label{sec:esoms}

In order to improve the understanding of the results discussed in \autoref{sec:comparison}, this section provides an overview (\autoref{sec:esoms1}) of energy system optimization models and general classification schemes (\autoref{sec:esoms2}), and develops a classification scheme for the present review (\autoref{sec:esoms3}).

\subsection{Overview} \label{sec:esoms1}
    \acp{ESOM} are built to retrieve a set of decisions on the operation and expansion of an energy supply system. These decisions compose a strategy which is supposed to be optimal with respect to a predefined objective. These kind of optimization models are distinguished from other types of energy system models in that their application constitutes mathematically optimized prescriptive analytics. 
    
    A systematic literature review on national energy system optimization modelling for decarbonization pathways by \citet{plazas2022national} classifies energy system models and lists MARKAL, IKARUS, OPERA, LUT Energy System Transition Model, TIMES, MESSAGE, OSeMOSYS, GENeSYS-MOD, TEMOA and EnergyScope as major models in the literature. Previous reviews of energy systems models are listed in the survey on energy systems modeling by \citet{pfenninger2014energy}. Furthermore, a review targeted specifically at open source energy model development by \citet{Groissbck2019AreOS} ranks them based on an evaluation of a weigthed degree of implementation of 81 proposed functions. The recent review of energy system optimization frameworks for model generation by \citet{Hoffmann2024ARO} contains a comparison of 63 energy system optimization frameworks according to the Open Energy Platform and the Open Energy Modelling Initiative. The latter, more recent study expands the list of major models to include Calliope, PyPSA, ETHOS.FINE and oemof, among others.

\subsection{General Classification Scheme} \label{sec:esoms2}
    We aim to identify a classification scheme tailored to linear \acp{ESOM}, which are a subset of the broader category of energy system models. According to \citet{beeck2001classification}, energy system models can be classified on nine dimensions: Purpose, Assumptions, Analytical Approach, Methodology, Mathematical Approach, Geographical Coverage, Sectoral Coverage, Time Horizon and Data Requirements. The \textbf{purpose} of a model is understood as the questions it addresses such as forecasting, exploration of the current system or assessments of different policies. \textbf{Assumptions} are distinguished to be about endogenous parameters of the model and exogenous ones that are supplied by the user. The \textbf{analytical approach} is the distinction between bottom-up models which build the system up from a detailed technological description and disaggregated data, and top-down models which describe the energy system from a macro-economic perspective with aggregated data and elasticities. The \textbf{methodology} describes the type of modelling applied such as econometric models, simulation or optimization models. The mathematical framework employed to build the model, such as linear programming, integer programming or differential equations are covered by the \textbf{mathematical approach}. The other dimensions describe spatial and temporal properties of the model and the data requirements classified into aggregated, disaggregated, quantitative and qualitative data. 
    
    Proceeding our exploration of classification schemes for energy system models, the publication by \citet{Mougouei2017EffectiveAT} offers a comprehensive overview of the existing literature on classification schemes. According to this study, the fundamental characteristics of energy system optimization models are: Analytical Approach, Mathematical Approach, Geographical Coverage, Sectoral Coverage and Time Horizon. A review on classifications of bottom-up models by \citet{Prina2020ClassificationAC} proposes the mathematical approach as well as coverage and resolution of different dimensions, together with the information on the modelling technique and the type of decisions to be optimized. A broad review on energy system models by \citet{Klemm2021ModelingAO} classifies in two directions: One direction takes the analytical and mathematical approach into consideration, among additional categories, to account for models that are not related to optimization and characteristics related to their usability and purpose. The other direction takes the model's technological granularity into consideration, which includes spatial, temporal, sectoral and economic coverage as well as resolution. The 2014 founded Open Energy Modelling Initiative \cite{pfenninger2018opening} characterizes every model in their wiki according to the following dimensions: Model class that represents the mathematical approach, covered sectors, technologies included, decisions which are either of type dispatch or investment, scope of regions, geographic resolution, time resolution, network coverage and type of uncertainty modeling. 

\subsection{Review Classification Scheme} \label{sec:esoms3}
    In the previous survey on general classification schemes, we recognized the following pattern: The model's approach, both analytical and mathematical, are two basic dimensions to take into consideration (see \autoref{tab:classification}). Beyond that, all schemes take a subset of the spatio-temporal and the techno-economic scope and resolution into account. Finally, the type of decisions are of interest, especially due to their relation to the model's mathematical structure. 

    Both analytical and mathematical approach fall into a nominal scale, therefore yielding corresponding classes (\autoref{tab:classification}). Proceeding on the mathematical approach, we only take linear models into account, therefore LP, MILP and ILP models need to be represented, where the parameters can be deterministic and stochastic. Based on the model's mathematical properties, we added a classification of decisions, which can include expansion or scheduling-type decisions, or both. Expansion decisions deal with investments for capacity expansion planning, transmission expansion planning or generation expansion planning (often over several investment periods leading to transformation pathways) while scheduling decisions cover economic dispatch, unit commitment schedules or optimal power flow dispatch. 
        
    Based on the general classification schemes, we can also operationalize the scope of the model's spatial, temporal, technological and economic dimensions. Most schemes define classes such as "low", "medium" and "high", which need to be well-defined in order to make models comparable. Their meaning depends on whether the definition uses a relative or an absolute measure. In order to simplify that definition, we employ a binary classification, only reporting if the model contains dimensions of these types. A single-node model would not contain spatial dimensions and a single-sector model has no further economic dimensions. On the contrary, all of the \acp{ESOM} contain technological dimensions to account for the different components modelled. All symbols for the classifications are given in \autoref{tab:classification} and used when describing the results in \autoref{sec:comparison}.

    \begin{table*}
        \caption{Linear \acs{ESOM} Classification Scheme for this review.}
        \label{tab:classification}
        \tiny
        \begin{tabular}{|m{0.25\linewidth}|m{0.33\linewidth}|m{0.33\linewidth}|}
            \hline
             \cellcolor{\tableHeaderColor}\textbf{Analytical Approach} & \cellcolor{\tableHeaderColor}\textbf{Mathematical Approach} &
             \cellcolor{\tableHeaderColor}\textbf{Scope} \\ \hline
             \multirow{2}{=}{\begin{itemize}[leftmargin=*,noitemsep,topsep=0pt]
                \item \textbf{(TD)}: Top-Down Approach
                \item \textbf{(BU)}: Bottom-Up Approach
             \end{itemize}} & 
             \multirow{2}{=}{\begin{itemize}[leftmargin=*,noitemsep,topsep=0pt]
                \item \textbf{(C)}: Continuous variables present
                \item \textbf{(I)}: Integer variables present
                \item \textbf{(S)}: Stochastic parameters
             \end{itemize}} &
             \begin{itemize}[leftmargin=*,noitemsep,topsep=0pt]
                \item \textbf{(s)}: spatial dimensions
                \item \textbf{(t)}: temporal dimensions
                \item \textbf{(e)}: economic dimensions
             \end{itemize} \\ \cline{3-3}
             &
             &
             \begin{itemize}[leftmargin=*,noitemsep,topsep=0pt]
                \item \textbf{(SH)}: Scheduling type decisions
                \item \textbf{(EX)}: Expansion type decisions
             \end{itemize} \\ \hline
        \end{tabular}
    \end{table*}

%% file: sections/05-parallelization.tex
\section{Parallel Decomposition}
\label{sec:parallelization}

In order to parallelize the search for the optimum, a problem needs to be decomposed into independent pieces which can be processed by different processes. This can be achieved through parallelized model decomposition or parallelization of the solver's units. 

The parallelization of the algorithm's functional units takes place within the subroutines of the solver and usually encompasses parallelized Cholesky factorization or parallelized versions of the numeric procedures that find solutions to the given linear equations, usually taking advantage of data parallelism, e.g. \cite{rehfeldt2022massively} or using \ac{SIMD} processing, e.g. \cite{hafsteinsson1994solving}. If the algorithm can be decomposed into individual tasks such that each one processes a different part of the input, task parallelism is obtained.

Decomposing the domain model into sub-models leads to task parallelism. The model is decomposed according to a structure that can be either identified in its algebraic formulation or by investigation of the non-zero patterns in its generated constraint matrix. Given that substructure, a corresponding decomposition method can be selected \cite{conejo2006decomposition, ConstanteFlores2025OptimizationVR}.

\begin{figure}
    \centering
    \begin{subfigure}[t]{.5\textwidth}
        \centering
        \caption{}
        \label{fig:h-struct}
        \includegraphics[width=\linewidth]{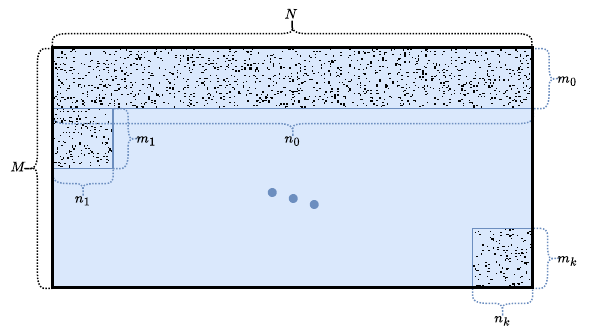}
    \end{subfigure}\hfill
    \begin{subfigure}[t]{.5\textwidth}
        \centering
        \caption{}
        \label{fig:v-struct}        
        \includegraphics[width=\linewidth]{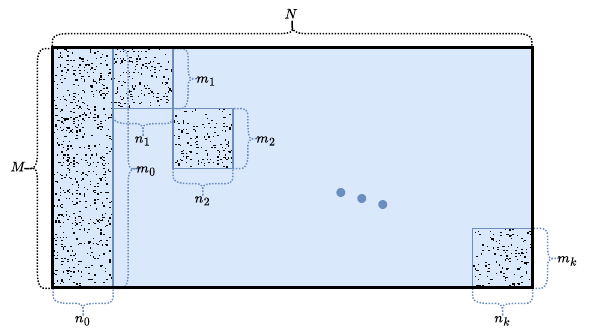}
    \end{subfigure}\hfill
    \begin{subfigure}[t]{.5\textwidth}
        \centering
        \caption{}
        \label{fig:arrow-struct}
        \includegraphics[width=\linewidth]{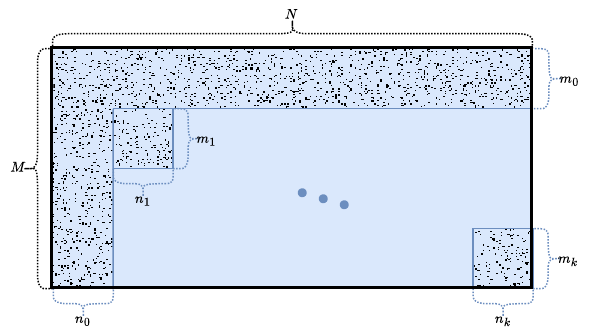}
    \end{subfigure}\hfill
    \begin{subfigure}[t]{.5\textwidth}
        \centering
        \caption{}
        \label{fig:stair-struct}
        \includegraphics[width=\linewidth]{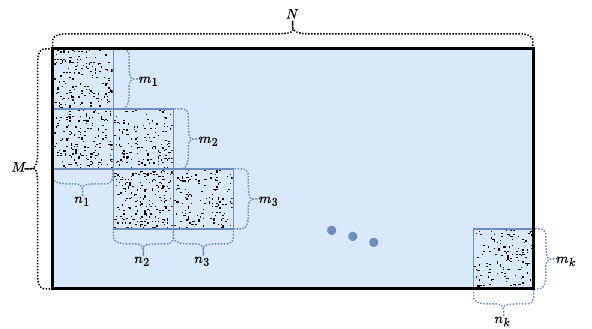}
    \end{subfigure}
    \colorbox{\legendColor}{
        \tiny
        \begin{tabular*}{.9\linewidth}[t]{r@{:\enspace}l}
             $(N,M)$ & (columns, rows) of the matrix \\
             $(n_0,m_0)$ & (columns, rows) of the master-block \\
             $(n_i,m_i)$ & (columns, rows) of $i^{th}$ block
        \end{tabular*}}
    \caption[]{ Different patterns of coefficient matrices.  
        \begin{tabular}[t]{r@{:\enspace}l}
            (\subref{fig:h-struct}) & Horizontally Bordered Block-Diagonal Matrix. \\
            (\subref{fig:v-struct}) & Vertically Bordered Block-Diagonal Matrix. \\
            (\subref{fig:arrow-struct}) & Arrowhead Matrix. \\
            (\subref{fig:stair-struct}) & Staircase Matrix.
        \end{tabular}}
    \label{fig:sub-struct}
\end{figure}

This section outlines decomposition methods for linear models, which can be expressed as
\begin{equation}
    min_{\bm{x}}\left(~\bm{c}\cdot\bm{x}~\mid~\bm{x}\in P\cap (\mathbb{R}^p\times \mathbb{Z}^q)~\right)
\label{eq:standardMILP}
\end{equation}
for polyhedron 
\begin{equation}
    P=\left\{~\bm{x}\in\mathbb{R}^N~\mid~\textstyle\sum_{j}A_{ij}x_j\leq b_i~\right\}
\label{eq:polyhedron}
\end{equation}
with coefficient matrix $\bm{A}$, bounding vector $\bm{b}$, cost vector $\bm{c}$ and decision vector $\bm{x}$. For the decision vector, $p$ entries are from $\mathbb{R}$ and $q$ entries are restricted to $\mathbb{Z}$ and $p+q=N$. If $q=0$, we deal with a convex LP, otherwise with a MILP that has a convex LP relaxation when integrality is dropped by relaxing $\mathbb{Z}^q$ into $\mathbb{R}^q$. The coefficient matrix may have a specific pattern created by its non-zero entries as shown in \autoref{fig:sub-struct} and \autoref{fig:opt-struct}. The model may yield such a substructure, either by algebraic construction or by permutation of rows and columns of the constraint matrix, which is possible due to the commutativity of the linear forms. The permutations are to be done in the tableau form as to not lose the association to the bounding vector and the cost vector. A corresponding decomposition method can be applied if a substructure is present. For a singly bordered substructure as shown in \autoref{fig:h-struct} and \autoref{fig:v-struct}, the coupling master-block is denoted as the $m_0\times n_0$ matrix $\bm{A}_0$.

    \subsection{Optimally Decomposable Substructure}
\begin{figure}
            \centering
            \begin{minipage}[l]{.49\textwidth}
                \centering
                \includegraphics[width=\textwidth]{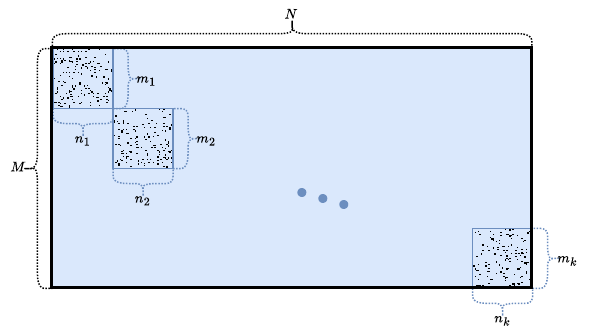}
                \colorbox{\legendColor}{
                    \tiny
                    \begin{tabular*}{.8\textwidth}[t]{r@{:\enspace}l}
                         $(N,M)$ & (columns, rows) of the matrix \\
                         $(n_i,m_i)$ & (columns, rows) of $i^{th}$ block
                    \end{tabular*}}
                \caption{Block-Diagonal Coefficient Matrix.}
                \label{fig:opt-struct}
            \end{minipage}    
            \begin{minipage}[r]{.5\textwidth}
                \centering
                \begin{minipage}[l]{.45\textwidth}
                    \centering
                    \begin{framed}
                        \begin{subfigure}[t]{.5\textwidth}
                            \centering
                            \includegraphics[width=\textwidth]{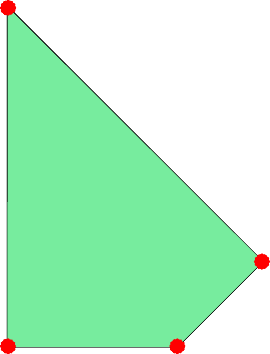}
                            \caption{}
                            \label{fig:opt-facet}
                        \end{subfigure}
                    \end{framed}
                    \begin{framed}
                        \begin{subfigure}[b]{.5\textwidth}
                            \centering
                            \includegraphics[width=\textwidth]{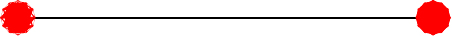}
                            \caption{}
                            \label{fig:opt-line}
                        \end{subfigure}
                    \end{framed}
                \end{minipage}
                \begin{subfigure}[r]{.47\textwidth}
                    \begin{framed}
                        \includegraphics[width=\textwidth]{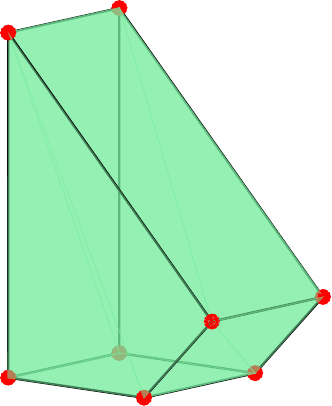}
                        \caption{}
                        \label{fig:opt-product}
                    \end{framed}
                \end{subfigure}
                \caption{Factorization (\subref{fig:opt-facet}, \subref{fig:opt-line}) of a separable polytope (\subref{fig:opt-product}).}
                \label{fig:opt-polytope}
            \end{minipage}
        \end{figure}

        All non-zero entries form blocks within the constraint matrix such that no two blocks overlap on any axis, as shown in \autoref{fig:opt-struct}. The dimensions $(m_i,n_i)$ of different blocks might significantly vary but each block represents an individual optimization problem of its own.
        An example of a separable coefficient matrix describing a polytope in three dimensions is given by \autoref{eq:opt-matrix}, together with the corresponding sub-blocks. 

        \begin{equation}
            [\bm{A}|\bm{b}] = \left[
                \begin{array}{ccc|c}
                    \cellcolor{\matrixBlockColor}{1} & \cellcolor{\matrixBlockColor}{1} & 0 & 1 \\
                    \cellcolor{\matrixBlockColor}{1} & \cellcolor{\matrixBlockColor}{-1} & 0 & \frac{1}{2} \\
                    0 & 0 & \cellcolor{\matrixBlockColor}{1} & \frac{1}{2}
                \end{array}
            \right],
            \bm{A}_1=\begin{bmatrix}1&1\\1&-1\end{bmatrix},
            \bm{A}_2=\begin{bmatrix}1\end{bmatrix}
            \label{eq:opt-matrix}
        \end{equation}

        ~

        Each block, together with the bounding vector, defines a polyhedron in a subspace of the full polyhedrons's space. The subspace for block $\bm{A}_1$ is two-dimensional as the block has two columns, while $\bm{A}_2$ describes a line in a one-dimensional subspace.
        The full polyhedron is the product of these two orthogonal factors, illustrated via the polyhedral geometry research software \textit{polymake} \cite{gawrilow2000polymake} in \autoref{fig:opt-polytope}. All of the sub-problems expressed by each independent block can be solved in parallel. Projecting the full polyhedron's optimum onto the orthogonal factors, these projections are going to be equal to the factors' optima. Therefore, the full polyhedron's optimum is simply reconstructed by vector addition of the factors' optimal points.

    \subsection{Constraint-Coupled Decomposable Substructure}
        If there exist constraints which couple variables of more than one block within the constraint matrix, those constraints, when grouped together, will form a block along one of the horizontal borders of the constraint matrix, shown in \autoref{fig:h-struct}. Here, all constraints which belong to the horizontal master-block couple at least two different sub-blocks, rendering them interdependent. This also means that those blocks are not orthogonal anymore. Projecting the polyhedron's optimum on such a block's dimensions could render a point outside the original polyhedron's region. Therefore, the superposition of the sub-blocks' solutions is not guaranteed to be feasible for the full polyhedron. In order to incorporate the coupling constraints while solving sub-blocks independently, we can employ the Dantzig-Wolfe decomposition or the Lagrange Relaxation technique. \\

        \noindent\textit{Dantzig-Wolfe Decomposition} \cite{dantzig1960decomposition}: The idea of the Dantzig-Wolfe Decomposition is a reformulation of the structured model such that we change the polyhedron's $\mathcal{H}$-representation into its $\mathcal{V}$-representation. This reformulation turns the master-block into a weighting problem, as it is reformulated into the convex combination of the all vertices of $P$ and, if $P$ is unbounded in some directions, the conical combination of those directions. The master block optimizes the weights of those combinations such that the resulting point is optimal in terms of the objective vector. As all vertices cannot be enumerated in practice, we start with a reduced master problem that contains only a few ones. This reduced $\mathcal{V}$-representation of the polytope can be expressed as a convex combination $\{\bm{x}\in\mathbb{R}^N~\mid~\mathbf{V}\bm{\lambda}=\bm{x}~\land~ \Vert\bm{\lambda}\Vert_1=1~\land~\lambda_i\geq 0\}$ for some matrix $\mathbf{V}$ that contains a subset of the polyhedron's vertices as its columns. This is the weighting problem. The sub-blocks on the contrary provide these vertices to the master-block. In order to get vertices which are driving us to the global optimum, we need to adjust each sub-block's objective vector such that, given the current convex combination in the master-block still is not globally optimal, the next set of vertices is improving our master-solution when included into its convex combination. The sub-blocks are referred to as the pricing problems. Basically, the reduced master problem corresponds to a partial polyhedron which is iteratively expanded into the direction of the optimal vertex. This method, referred to as delayed column generation, is supposed to terminate before enumerating all vertices and thus saving both runtime and memory. All sub-problems can be requested in parallel for every iteration and need to be synchronized at the beginning of each new iteration. 
        
        In the case of a MILP, these ideas have to be incorporated into the Branch-And-Bound procedure leading to the Branch-And-Price method, comprehensively explained in \cite{G-2024-36}. \\
        
        \noindent\textit{Lagrange Relaxation} \cite{Geoffrion1972LagrangeanRA}: In contrast to the Dantzig-Wolfe decomposition, Lagrange relaxation expands the feasible region of the polyhedron by removing the master-block constraints. This expansion is controlled by Lagrange multipliers $\bm{u} \in \mathbb{R}^{m_0}$, also known as dual variables, which penalize violations of the relaxed constraints. Using these multipliers, we can define the Lagrangian function as follows:
        
        \begin{equation}
            \mathcal{L}(\bm{x}, \bm{u}) = \sum_{i=1}^N c_i x_i + \sum_{i=1}^{m_0} u_i \cdot \left[\sum_{j=1}^{N+m_0} A'_{ij} x'_j - b_i\right]
        \end{equation}
        with $\bm{A}'$ and $\bm{x}'$ as the coefficient matrix and decision vector in slack form.
        From this, the dual function is derived, which maps each vector of dual variables to the optimal value of the Lagrangian function:
        \begin{equation}
            d(\bm{u}) := \underset{\bm{x}}{\text{inf}} \; \mathcal{L}(\bm{x},\bm{u})
        \end{equation}
        The resulting dual function is always concave as it is the lower envelope of the family of Lagrangians and in general, it is continuous but not differentiable \cite{boyd2004convex, bagirov2020numerical}. The dual function provides a lower bound on the objective value of the primal problem. Therefore, the goal becomes to find the best possible lower bound, which leads to the dual problem: $\max_{\bm{u}}(d(\bm{u}) \;|\; \bm{u} \geq \bm{0})$.
        Due to the nature of the primal constraints, the dual variables are typically required to satisfy: $\bm{u} \geq \bm{0}$.
        Because the dual function is not differentiable in general, solving the dual problem requires methods from nonsmooth optimization. In this context, only a subgradient can be computed. The subgradient method, which generalizes gradient descent, has the key challenge that a subgradient at the optimum does not necessarily vanish, which makes it difficult to verify optimality \cite{bagirov2020numerical}. Its main advantage lies in its simplicity.
        
        Another common approach is the bundle method. In this method, after each evaluation of the dual function, first-order information is stored in a so-called bundle. A piecewise linear approximation of the dual function is constructed from this bundle, and the optimal value of this approximation is sought. This method often leads to better convergence properties than subgradient approaches \cite{bagirov2020numerical}.
        
        A further technique is to improve the convergence by stabilizing the iterates and penalizing large constraint violations more strongly through a regularization term. Often, that term corresponds to a proximal operator, a mapping that generalizes projection by balancing objective minimization with proximity to a reference point. This leads to an augmented Lagrangian upon which the Alternating Direction Method of Multipliers builds up to solve structured problems \cite{beck2017first}. 
    
    \subsection{Variable-Coupled Decomposable Substructure}
        If there exist variables which couple constraints of more than one block within the constraint matrix, those variables, when grouped together, will form a block along one of the vertical borders of the coefficient  matrix as \autoref{fig:v-struct} illustrates.
        
        Here, all variables which belong to the vertical master-block couple at least two different sub-blocks, rendering them interdependent. This again means that the blocks are not orthogonal anymore. Projecting the polyhedron's optimum on such a block's dimensions could render a point outside the original polyhedron's region. Therefore, the superposition of the sub-blocks' solutions is not guaranteed to be feasible for the full polyhedron. In order to incorporate the coupling variables while solving sub-blocks independently, we can employ the Benders Decomposition or the Variable Splitting technique.

        ~

        \noindent\textit{Benders Decomposition} \cite{bnnobrs1962partitioning}: The Benders decomposition splits the original problem into a master problem and one or more sub-problems. The original problem is projected to the subspace that is defined by its coupling variables \cite{rahmaniani2017benders}. The resulting formulation is dualized, resulting in an equivalent problem which only depends on the coupling variables. The other set of variables is replaced by associated cuts which represent the feasible space and the projected costs. However, too many constraints are added to solve the resulting problem directly \cite{martin2012large}.
        Therefore, the problem is relaxed by removing the cuts, resulting in the relaxed master problem. It represents the relaxed polyhedron and optimizes the coupling variables. The optimized variables from the master problem are fixed within the sub-problems, where the remaining variables are optimized. If the proposed fixed variables from the master problem result in feasible sub-problems, the sub-problems return optimality cuts to the master problem. Otherwise, the master problem receives feasibility cuts. These Benders cuts that are added to the master problem in each iteration approximate the objective function of the sub-problems from below. The sub-problems provide an upper bound to the objective function of the original problem since they are restricted by the fixed values of the coupling variables \cite{conejo2006decomposition}. The master problem and the sub-problems are iteratively solved until the convergence tolerance is reached.

        While the algorithm was initially formulated for mixed-integer linear programming models \cite{bnnobrs1962partitioning}, it has since been generalized, e.g. for nonlinear sub-problems \cite{geoffrion1972generalized}. Benders decomposition has been applied to a broad set of optimization models including  stochastic, bi-level and multi-stage problems \cite{rahmaniani2017benders}.

        ~
        
        \noindent\textit{Variable Splitting} \cite{guignard1987lagrangean}:
        Variable splitting is a reformulation technique used to facilitate the decomposition of variable coupled optimization problems, particularly in the context of Lagrangian relaxation . The core idea is to duplicate the coupling variables across sub-problems in order to separate the model. These duplicated variables are then linked via consensus constraints, effectively transforming the vertical master-block into a horizontal one.
    
    \subsection{Arrowhead Substructure}
        If there exist both variables and constraints that couple different blocks, the corresponding structure takes on the shape of an arrowhead (see \autoref{fig:arrow-struct}). In this case, we can use Variable Splitting to incorporate the vertical block into the horizontal one and employ a corresponding method for constraint-coupled decomposable substructures.        
    
    \subsection{Staircase Substructure}
       If the non-zero elements of the constraint matrix are only found on the diagonal blocks and adjacent off-diagonal blocks, the pattern resembles a staircase as shown in \autoref{fig:stair-struct}. This special structure can be exploited by compact basis methods, nested methods based on the Dantzig-Wolfe decomposition or a specifically adapted simplex method \cite{Fourer1982SolvingSL, Fourer1983SolvingSL}.

%% file: sections/06-implementations.tex
\section{Implementations}
\label{sec:comparison}

This section provides an overview of the existing literature on the parallelization of decomposition in energy system optimization problems (\autoref{sec:comparison1}), formulates recommendations for benchmarks of linear ESOMs (\autoref{sec:comparison2}), and provides an overview of the available software (\autoref{sec:comparison3}).

\subsection{Review of Parallelized Decomposition} \label{sec:comparison1}
    In this section we review the literature that has been considered eligible according to the process described in \autoref{sec:methods}. After $126$ records have been collected from two large database platforms, they have undergone a clearly set up screening process. The high amount of exclusions is due to the fact that a large amount of publications have not parallelized the decomposition method ($n=31$), or the term decomposition method was not referring to the exact methods described in \autoref{sec:parallelization}, being (meta-)heuristics, custom methods or decomposition in its broadest meaning ($n=27$). When using a framework to create the model, the configuration to set up the model is not necessarily reported. Also, meta-parameters to assess the size and complexity of the model were frequently excluded in the reporting, such as number of variables, number of constraints, model sparsity or integrality fraction. Considering the benchmarking methods, we have not encountered any publication that has sampled over several runs of a single solve in order to account for the performance variability of the underlying system. 
    
    The publications that have used a parallelized version of the decomposition method can be divided into two groups according to the type of benchmarking experiment that has been conducted. Either, the parallelized decomposition method has been compared to the system solved as a centralized program or the parallelized decomposition's scaling behaviour has been investigated by measuring its performance for different amounts of parallel computation processes. The first group is referred to as methodological benchmark and the second one as scaling benchmark. The methodological benchmark studies which employed a parallelized model decomposition method are shown in \autoref{tab:results-method-stratified}, with the type of model encoded in columns one to four according to the classification scheme in \autoref{tab:classification}. The full results of the data extraction from the included studies are listed in \ref{app:data}

    \tiny
    \begin{longtable}{
    R{0.12\linewidth}
    R{0.20\linewidth}
    R{0.10\linewidth}
    R{0.10\linewidth}
    R{0.17\linewidth}
    R{0.15\linewidth}}
        \caption{Methodological benchmarks on parallelized ESOM decomposition.}
        \label{tab:results-method-stratified}
        \\ \hline
        \cellcolor{\tableHeaderColor}\textbf{Type} &
        \cellcolor{\tableHeaderColor}\textbf{Name} &
        \cellcolor{\tableHeaderColor}\textbf{Size} &
        \cellcolor{\tableHeaderColor}\textbf{\mbox{Decomposition} \mbox{Dimensions}} &
        \cellcolor{\tableHeaderColor}\textbf{Decomposition Method} &
        \cellcolor{\tableHeaderColor}\textbf{Study}
        \csvreader[
            on column count error,
            head to column names,
            respect percent,
            respect underscore
        ]{results/DataExtraction-MethodBenchmark.csv}{}
        {\\ \hline \ESOMType & \ESOMName & \ESOMSize & \DecompositionDimensions & \DecompositionMethod & \Citation}
        \\ \hline
    \end{longtable}
    \normalsize
    
    Among the studies in \autoref{tab:results-method-stratified}, eleven reported a speedup that is larger than one \cite{Dvorkin2018,Liu2015,Fu2013,Zhang2024,Wu2010,Huang2017,Steven2024,Wales2024,Soares2022,Alhaider2018,Gke2024}. Two studies \cite{Paul2023,Gke2024} also observed speedup values lower than one. The reference solver often did not prove optimality either, and in instances reported by \citet{Wales2024}, improved solution accuracy was demonstrated. Furthermore, \citet{Wu2010} as well as \citet{Huang2017} could solve instances for which the reference solver's method did not converge for the given computational resources. 
    
    Among the methodological benchmarks which reported a quantitative quality metric for both the reference system and the system under test, no improvement in accuracy greater than $0.9\%$ for the optimality gap was observed \cite{Wales2024}, and the relative deviation from the known optimum increased by no more than $0.2\%$ \cite{Paul2023}. Those methodological benchmarks which reported a qualitative  quality metric for both the reference system and the system under test did not measure any degradation \cite{Wu2010,Huang2017,Gke2024}. These results indicate that parallelized decomposition methods for linear \acp{ESOM} do improve the optimization procedure's performance. However, most of the procedures have been tailored to improve beyond the textbook implementation of the corresponding decomposition method. Those improvements are particularly necessary if the degree of parallelization is low.
    
    Scaling benchmarks on parallelized decomposition methods have been performed by \citet{Gil2016}, \citet{Gong2019} and \citet{Sundarraj1995}. We want to highlight that the scaling benchmarks do not necessarily contain the performance metric for a single process run such that the estimation of the parameters from Gunther's law needs to be adapted by introducing another parameter to accommodate for the missing baseline performance \cite{gunther1997practical}. Given the small sample size in the studies, the nonlinear fit lacks robustness and may not be statistically meaningful, therefore we haven't performed a regression analysis for them. We also observed non-uniform sampling of the system's scaling curve where the independent variable was set to values which are powers of two instead of equidistant values. This would bias the fitted model's accuracy.
    
\subsection{Recommendations for Benchmark Guidelines} \label{sec:comparison2}
    Following the previously outlined observations on benchmarking methods, we provide a few recommendations for benchmarks on linear \acp{ESOM}. The state of reporting results of parallel computing experiments has been investigated by \citet{hoefler2015scientific}, who developed a set of twelve rules which help to maintain reproducibility and improve interpretability. We recommend to consider this set when evaluating the results of the experiment. Apart from these, we also recommend to include the following fundamental items which apply in the context of benchmarking different solving methods for linear \acp{ESOM}: \\
    
    \textbf{ReBeL-E}:
    \textbf{Re}commendations for \textbf{Be}nchmarks on \textbf{L}inear \textbf{E}SOMs
    \begin{enumerate}
        \item \textbf{Stats:} \texttt{Bring in stable measurements.}\\
        How many samples for each optimization run have been obtained. These are necessary to account for the performance variability of the whole system. Any further metric should be reported as a location parameter such as the average value and its error metric such as the standard deviation.
        \item \textbf{Model:} \texttt{Benchmark a standard model.}\\
        A publicly accessible and pre-configured model. This could be either an IEEE power system test case which can be found for example in PandaPower \cite{thurner2018pandapower} or a pre-configured energy system model such as PyPSA-Eur \cite{horsch2018pypsa}. If a model generator or framework needs to be used, its parametrization should be made fully accessible.
        \item \textbf{Size:} \texttt{How many dimensions?}\\
        The size of the model in terms of dimensionality, i.e. total number of variables. Additional values such as number of scenarios and regions are optional.
        \item \textbf{Complexity:} \texttt{Quantify the complexity.}\\
        The sparsity of the coefficient matrix, i.e. the total number of non-zero elements, and the total amount of constraints. Additionally, the absolute number of integer variables or their share as a fraction on the total number of variables should be included when dealing with \ac{MILP} models.        
        \item \textbf{Solver:} \texttt{Introduce us to the solver.}\\ 
        Which solver and which version has been used, preferably with its configuration.
        \item \textbf{Modelling Environment:} \texttt{What is your language?}\\
        Which algebraic modelling language or interface has been used to pass the model to the solver, as those can have a substantial impact due to automatic reformulations done by the AML transpiler.
        \item \textbf{Quality Metric:} \texttt{How good was it?}\\
        Which metric has been used to assess the quality of the optimization result, preferably the duality gap or the MIP gap. These can also include relative deviation of a known optimum or similar metrics, however, it should always be reported for every measurement for both the reference system and the system under test, even when it has timed out. Values are always reported with their error metric.
        \item \textbf{Performance Metric:} \texttt{Did it run or walk?}\\
        Which metric has been used to assess the method's performance, such as runtime or speedup for every measurement for both the reference system and the system under test, even if the system exceeded the available memory. Values are always reported with their error metric.
        \item \textbf{Reference Point:} \texttt{Do not drop the origin.}\\
        When performing a scaling experiment, the single-core performance and quality should always be included. When performing a methodological experiment, a reference system is used as a comparison.        
    \end{enumerate}

We recommend these items to be always included into any benchmark that assesses the performance of a proposed method on optimizing a linear \ac{ESOM}. Nevertheless, benchmarks are difficult to compare if not conducted on a standardized benchmark set, similar to the MIPLIB suite \cite{miplib2027}. Such a standard benchmark suite would be highly beneficial for the field of energy system optimization.

\subsection{Software for Parallelized Model Decomposition} \label{sec:comparison3}
    A variety of techniques have been developed to deal with large-scale optimization problems. Certain tools leverage the specific structure of models generated by algebraic modelling systems to apply decomposition techniques, enhancing the solver’s efficiency in handling large-scale problems. \autoref{tab:modelling-tools} outlines modeling tools, specifying how they use model-specific structures to facilitate decomposition.
    According to the core functionalities we identified, the tools are grouped into three categories, acknowledging that the boundaries between these categories are not always clear-cut.

\subsubsection*{Parallelization Facilitators}
    Several tools assist with the specific requirements of distributed computation. These tools are often built upon algebraic modeling languages. The GAMS Grid Facility \cite{GAMSGridFac} supports distributed computation for models written in the GAMS language. Similarly, \texttt{parAMPL} \cite{parAMPL} provides parallel execution capabilities for models formulated in AMPL. The tool \texttt{mpi-sppy} \cite{mpi-sppy}, a successor to \texttt{PySP}, focuses on parallel computations in stochastic programming models represented as scenario trees in Pyomo. Meanwhile \texttt{StructJuMP} \cite{Huchette2014ParallelAM} focuses on block structured two-stage stochastic optimization problems that are solved in parallel on distributed memory system. The modelling framework \texttt{StochasticPrograms.jl} \cite{Biel2019EfficientSP} employs parallelized solvers for stochstic programming problems such as L-Shaped solvers, Progressive Hedging Solvers and Quasi-gradient solvers. We also classify \texttt{disropt} \cite{farina2020disropt} in this category. It is a Python-based framework designed to establish consensus on the optimal solution of a distributed optimization problem, where the objective function and constraints are distributed among multiple agents. Communication is restricted to neighboring nodes, typically implemented using the Message Passing Interface \cite{clarke1994mpi}. The \texttt{Ubiquity Generator Framework} \cite{Shinano2017TheUG} allows for a given Branch-And-Bound based mixed-integer linear programming solver to be instantiated and distributed among parallel processes in a manycore system or on a cluster computing environment.

\subsubsection*{Modularized Branch-And-Price-And-Cut frameworks}
    The tools \texttt{BaPCod} \cite{sadykov2021bapcod}, \texttt{Coluna} \cite{Coluna}, \texttt{GCG} \cite{GCG}, and \texttt{DIP} \cite{DIP} provide modular frameworks for Branch-And-Price-And-Cut algorithms, supporting common decomposition techniques such as Benders Decomposition, Lagrangian decomposition schemes, and Dantzig-Wolfe reformulation. Given the many variation points in these algorithms, the frameworks offer default implementations for all essential components, while also allowing for user-defined extensions. \texttt{Coluna} is implemented in Julia and builds upon JuMP and the MathOptInterface. \texttt{GCG} is part of the SCIP optimization suite. It automatically finds permutations of the constraint matrix which create a block structure by solving a combinatorial optimization problem. \texttt{DIP} is distributed through COIN-OR and comes with limited documentation. \texttt{BaPCod} is available for academic use via email upon request.

\subsubsection*{Generic Abstractions for Decomposition-Based Algorithms}
    This final category encompasses tools with the greatest internal diversity. What unites them is each of them providing abstractions aimed at facilitating the implementation of novel decomposition based algorithms proposed by their developers. \texttt{DSP} \cite{DSP} supports stochastic programs by offering algorithmic scaffolding for both serial and parallel versions of Dantzig-Wolfe reformulation, Lagrangian decomposition schemes, Integer Benders Decomposition, or solving in extensive form using the underlying solver. The framework \texttt{urbs-DecEnSys} \cite{urbsDecEnSys} targets energy system models involving time series. While it does not support integer variables, it facilitates parallel sub-problem solutions for linear problems such as capacity expansion and unit commitment. These decompositions can be based on time, spatial regions, Benders Decomposition, or Stochastic Dual Dynamic Programming.

    \texttt{Plasmo.jl} \cite{Jalving2022} is a Julia package designed to automate the identification of promising decomposition schemes. Its central concept is to represent optimization problems as hypergraphs, where nodes correspond to variables and hyperedges to constraints involving those variables. Hypergraph algorithms on the enhanced algebraic model are then used to detect decomposition opportunities. This approach enables hybrid decompositions across time and space, which may not be readily apparent to modelers. The companion package \texttt{PlasmoAlgorithms} \cite{cole2025} supports the implementation of decomposition strategies such as Benders Decomposition and Dual Dynamic Programming, enabling the evaluation and exploitation of these decomposition structures.

    Finally, \texttt{SMS++} \cite{SMSpp} offers the most flexible yet technically demanding infrastructure for constructing custom decomposition algorithms. Its developers address the gap between formulating mathematical models generically and a block-structured way that facilitates decomposition algorithms: \texttt{SMS++} supports decomposition implementations via reusable and nestable abstractions applicable across a wide range of decomposition techniques. The most general of these abstractions is the \textit{Block}, an abstract base class for representing a self-contained part of the model. In addition, the \textit{Solver} abstraction can represent either an off-the-shelf solver or a custom algorithm designed to exploit specific structures in a \textit{Block}, potentially nesting them.

    \tiny
    \begin{longtable}{
    R{0.11\linewidth}
    R{0.08\linewidth}
    R{0.22\linewidth}
    R{0.08\linewidth}
    R{0.15\linewidth}
    R{0.06\linewidth}
    R{0.11\linewidth}}
        \caption{Software for Parallel Decomposition. \textit{Abbreviations: \acf{BB}, \acf{BPC}, \acf{LP}, \acf{MILP}, \acf{MIP}, \acf{MPI}, \acf{SP}}}
        \label{tab:modelling-tools}
        \\ \hline
        \cellcolor{\tableHeaderColor} \textbf{Name} & 
        \cellcolor{\tableHeaderColor} \textbf{Language} & 
        \cellcolor{\tableHeaderColor} \textbf{Core Concepts} & 
        \cellcolor{\tableHeaderColor} \textbf{License} & 
        \cellcolor{\tableHeaderColor} \textbf{Parallelization} & 
        \cellcolor{\tableHeaderColor} \textbf{Problem Types} & 
        \cellcolor{\tableHeaderColor} \textbf{Last Update as of May 2025}
        \\ \hline
        \textbf{BaPCod} & C++ & Highly customizable \acs{BPC} scheme & Academic EULA & \acs{MIP} solver parallelization through Multi-Threading & \acs{MIP} & Aug 2024\\
        \hline
        \textbf{Coluna} & Julia & Automatic \acs{BPC} for block-structured \acsp{MIP} &  MPL 2.0 & Can solve subproblems in parallel& \acs{MIP} & Feb 2024\\
        \hline
        \textbf{DIP} & C++, Python interface  & Provides customizable algorithmic implementation details for \acs{BPC} and related decomposition-based algorithms&  EPL 1.0 & Subproblems and \acs{BB} tree via the Abstract Library for Parallel Search \cite{xu2005alps} & \acs{MIP} & Jan 2020\\
        \hline
        \textbf{disropt} & Python & Distributed optimization agents eventually reach consensus about optimal solution & GPL 3.0 & Unlimited amount of agents& Any  & Jun 2021\\
        \hline
        \textbf{DSP} & C++ & Serial and Parallel implementation of Dantzig-Wolfe, Lagrange, and Benders Decompositions with CPLEX, Gurobi, or SCIP underlying& 3-clause BSD & \acs{MPI} & \acs{MIP} & Jun 2023\\
        \hline
        \textbf{urbs-DecEnSys} & Python & Provides energy-specific abstractions for modeling distributed energy systems of any scale with time series data &  GPL 3.0 & Only in underlying \acs{LP} solver through pyomo & \acs{LP} & Jul 2019\\
        \hline
        \textbf{GAMS Grid Facility} & GAMS & Facilitates asynchronous submission and collection of GAMS-model solution tasks on HPC Grids and multi-core systems & Academic and commercial EULA & Native & Any & Mar 2025\\
        \hline
        \textbf{GCG} & C++, other interfaces & Automatic structure detection with Dantzig-Wolfe reformulation or Benders, modularized \acs{BPC}& GNU LGPL & Can solve pricing in parallel&\acs{MIP} & Apr 2025\\
        \hline
        \textbf{mpi-sppy} & Python & Pyomo extension that supports scenario-discretization of multi-stage stochastic programs &  3-clause BSD &  \acs{MPI} & \acs{SP} & May 2021\\
        \hline
        \textbf{parAMPL} & Python & Facilitates asynchronous submission and collection of AMPL-model solution tasks on HPC Grids and multi-core systems & 2-clause EULA & Native & Any & Sep 2019\\
        \hline
        \textbf{Plasmo.jl} & Julia & Structure identification via hypergraph description of optimization problem &  MPL 2.0 &  Interface to PIPS-NLP, can solve subproblems in parallel& Any & Nov 2024\\
        \hline
        \textbf{SMS++} & C++ & Providing software abstractions more suitable to decompositions than general algebraic modeling abstractions &  GNU LGPL v3 & Specific to each \texttt{Solver} configured for a \texttt{Block} & Any & May 2025\\
        \hline
        \textbf{StructJuMP} & Julia & Define blocks and linking variables explicitly to distribute the decomposed model in parallel & MIT Expat License & \acs{MPI} & Any & Nov 2023\\
        \hline
        \textbf{Stochastic Programs.jl} & Julia & Different abstractions for the core elements of \ac{SP}s to form blocks & MIT License & Can solve in parallel & \acs{SP} & Sep 2022\\
        \hline
        \textbf{Ubiquity Generator Framework} & C++ &  Parallelization of \ac{BB} solvers for distributed or shared memory systems &  GNU LGPL v3 & Automatic coordination of parallel solver instances  & \ac{MILP} & Nov 2024\\
        \hline
    \end{longtable}
    \normalsize

\subsubsection*{Software for Parallel Decomposed Model Solving}
    Given an already block-structured coefficient matrix of a linear optimization model, some solvers can directly exploit this property within linear algebra subroutines. They parallelize the solver's functional parts in order to speed up the overall computation. 
    One solver that takes advantage of the singly bordered block-structure is \texttt{PIPS-IPM} \cite{lubin2011scalable}. Assuming this structure, the computationally expensive step of solving a large linear system in interior point methods can be accelerated by parallelizing the numerical solution of a system of linear equations. This system arises from the Karush-Kuhn-Tucker conditions and the block-structure of the coefficient matrix propagates into this system. This method has been extended to arrowhead structures resulting in \texttt{PIPS-IPM++} by \citet{rehfeldt2022massively}. Another structure exploiting solver is \texttt{DuaLip} \cite{gupta2023practical,basu2020eclipse}, which solves a perturbed version of the linear program via gradient-based algorithms on its smooth dual. It also assumes a block diagonal structure in the coefficient matrix of the problem. Given integrality conditions, a highly distributed Branch-And-Bound solution process is employed by \texttt{ParaXpress} \cite{ParaXpress} as well as \texttt{ParaSCIP} \cite{ParaScip} and \texttt{FiberSCIP} \cite{FiberSCIP} which use the \texttt{Ubiquity Generator Framework} \cite{Shinano2017TheUG}. Recent developments on the primal-dual hybrid gradient method show improved performance for \acs{LP} problems \cite{Applegate2021, applegate2025pdlp} on GPUs and have been included into several software systems, highlighted in an overview on first-order methods parallelized on GPU devices by \citet{lu2025overview}.

%% file: sections/07-conclustion.tex
\section{Summary and Conclusion}
\label{sec:conclusion}
As the complexity and dimensionality of energy system optimization models continue to grow in response to increasing renewable integration, sectoral coupling, and spatio-temporal granularity, traditional solution techniques reach computational limits. This systematic review demonstrates that decomposition methods which exploit identifiable block structures are suitable for scaling linear energy system optimization with high-performance computing methods. This review processed 126 publications in total, out of which 15 publications matched the inclusion criteria yielding 79 benchmark instances. 
We found that no single decomposition method universally dominates and the suitability of a technique depends on the structural characteristics of the model.
It is also important to find a standard benchmark set in order to assess the performance of new methods for their general use in energy system optimization. This review reveals critical gaps in reproducibility, standardization, and benchmarking rigor. Therefore, we strongly advocate for the adoption of comprehensive and transparent evaluation protocols. Recommendations for conducting such studies have been developed in this publication and could benefit future benchmark studies. Also, the establishment of publicly accessible benchmark suites like MIPLIB, tailored for linear ESOMs, are necessary. Current efforts to implement such a benchmark suite are done by \citet{OpenEnergyBenchmark2025}.
Finally, while software ecosystems supporting decomposition and parallel solving have matured significantly, particularly frameworks like \texttt{UG}, \texttt{Plasmo.jl} and \texttt{SMS++}, there is currently no automated way to detect structures in a plain formulation without additional data on the structure yielded by the modeller, aside from \texttt{GCG} which needs to solve a computationally expensive optimization problem to find these structures.

A central priority for future studies should be the establishment of a standardized benchmark suite for linear energy system optimization models that reflects the diversity of energy system structures and scales. Such a repository should be accompanied by a rigorous, community-agreed protocol for reporting benchmark results. On the computational side, specialized linear algebra routines, tailored to structured problems, can significantly accelerate subproblem computations. These routines can be integrated into lightweight, open-source solvers such as TulipJL \cite{Tulip.jl}. Also, first order methods are gaining traction for their amenability to large-scale, highly parallel environments. Finally, hybrid computing architectures present a promising path for future research. The use of accelerator units (e.g. GPUs, FPGAs, TPUs) alongside general-purpose CPUs within distributed systems offers an opportunity to exploit parallelism at multiple levels: across submodels, within solver routines and within linear algebra kernels. Realizing this potential will require both algorithmic adaptation and efficient methods for their orchestration.